\documentclass[11pt,notitlepage,twoside,a4paper]{amsart}
\usepackage{mathrsfs}

\usepackage{amsmath,amssymb,enumerate}
\usepackage{epsfig,fancyhdr,color}
\usepackage{epstopdf}
\usepackage{amssymb}
\usepackage{amsmath,amsthm}
\usepackage{latexsym}
\usepackage{amscd}
\usepackage{psfrag}
\usepackage{graphicx}
\usepackage{epsf}
\usepackage[latin1]{inputenc}
\usepackage[all]{xy}
\usepackage{tikz}
\usepackage{prettyref}
\usepackage{subfigure}
\usepackage{float}
\usepackage{color}
\usepackage{array}
\usepackage{cite}
\usepackage[totalwidth=17cm,totalheight=24cm]{geometry}

\DeclareMathOperator{\arcsinh}{arcsinh}
\DeclareMathOperator{\arccosh}{arccosh}

\DeclareMathOperator{\trace}{tr}
\DeclareMathOperator{\inj}{inj}
\DeclareMathOperator{\area}{Area}

\newtheorem{theorem}{\rm\bf Theorem}[section]
\newtheorem{proposition}[theorem]{\rm\bf Proposition}
\newtheorem{lemma}[theorem]{\rm\bf Lemma}
\newtheorem{corollary}[theorem]{\rm\bf Corollary}
\newtheorem{definition}[theorem]{\rm\bf Definition}
\newtheorem{remark}[theorem]{\rm\bf Remark}

\newcommand{\N}{{\mathbb N}}
\newcommand{\R}{{\mathbb R}}

\newcommand{\cE}{{\mathcal E}}
\newcommand{\cT}{{\mathcal T}}
\newcommand{\cML}{{\mathcal{ML}}}
\newcommand{\supp}{{\mbox{supp}}}

\def\interieur#1{\mathord{\mathop{\kern 0pt #1}\limits^\circ}}

\title[Spacetimes with particles]{Constant Gauss curvature foliations of AdS spacetimes with particles}
\author{Qiyu Chen}
\address{Qiyu Chen:  
School of Mathematics, Sun Yat-Sen University,
510275, Guangzhou, P. R. China;
Mathematics Research Unit, BL G,
University of Luxembourg, Campus Kirchberg,
6, rue Richard Coudenhove-Kalergi,
L-1359 Luxembourg, Luxembourg.}
\email{chenqy0121@gmail.com}

\author{Jean-Marc Schlenker}
\address{Jean-Marc Schlenker:
Mathematics Research Unit, BL G,
University of Luxembourg, Campus Kirchberg,
6, rue Richard Coudenhove-Kalergi,
L-1359 Luxembourg, Luxembourg}
\email{jean-marc.schlenker@uni.lu}

\date{v1, \today}

\begin{document}

    \begin{abstract}
        We prove that for any convex globally hyperbolic maximal (GHM) anti-de Sitter (AdS) 3-dimensional space-time $N$ with particles (cone singularities of angles less than $\pi$ along time-like curves), the complement of the convex core in $N$ admits a unique foliation by constant Gauss curvature surfaces. This extends, and provides a new proof of, a result of \cite{BBZ2}. We also describe a parametrization of the space of convex GHM AdS metrics on a given manifold, with particles of given angles, by the product of two copies of the Teichm\"uller space of hyperbolic metrics with cone singularities of fixed angles. Finally, we use the results on $K$-surfaces to extend to hyperbolic surfaces with cone singularities of angles less than $\pi$ a number of results concerning landslides, which are smoother analogs of earthquakes sharing some of their key properties.


\thanks{Q. Chen is partially supported the by International Program for Ph.D. Candidates, Sun Yat-Sen University.}
	\end{abstract}

	\maketitle
	
    \section{Introduction}

    Let $\theta=(\theta_{1},...,\theta_{n_0})\in (0,\pi)^{n_0}$. In this paper we consider an oriented closed surface $\Sigma$ of genus $g$ with $n_0$ marked points $p_{1},...,p_{n_0}$ and suppose that
    \begin{equation*}
        2\pi(2-2g)+\sum_{i=1}^{n_0}(\theta_{i}-2\pi)<0.
    \end{equation*}

    This ensures that $\Sigma$ can be equipped with a hyperbolic metric with cone singularities of angles $\theta_{i}$ at the marked points $p_{i}$ for $i=1,...,n_0$ (see e.g. \cite{Troyanov}). Denote by $\mathcal{T}_{\Sigma,\theta}$ the Teichm\"uller space of hyperbolic metrics on $\Sigma$ with fixed cone angles, which is the space of hyperbolic metrics on $\Sigma$ with cone singularities of angle $\theta_{i}$ at $p_{i}$, considered up to isotopies fixing each marked point (see more precisely Section 2.1).

    \subsection{AdS spacetimes with particles}
        We are interested in 3-dimensional manifolds endowed with an AdS structure, that is, a geometric structure locally modeled on $AdS_3$, a complete 3-dimensional spacetime of constant curvature $-1$.
        Such AdS manifolds are also called AdS spacetimes, since they occur naturally in connection to gravitation. An AdS spacetime is {\em globally hyperbolic compact} (GHC) if it contains a closed Cauchy surface, and it is {\em globally hyperbolic compact maximal} (GHM) if in addition any isometric embedding into a globally hyperbolic compact spacetime of the same dimension is an isometry. GHM AdS spacetimes have been shown by G. Mess \cite{mess,mess-notes} to present remarkable analogies with quasifuchsian hyperbolic manifolds. Here we are particularly interested in AdS spacetimes with {\em particles}, that is, cone singularities of angles less than $\pi$ along time-like lines, as in \cite{BS}. Cone singularities of this type are used in the physics literature to model point particles in 3d gravity, see e.g. \cite{thooft1,thooft2}. (More details on AdS spacetimes with particles can be found in Section 2.2.)

        We say that a GHM AdS spacetime with particles is {\em convex} if it contains a convex Cauchy surface. Convex GHM AdS spacetimes with particles contain a smallest non-empty convex subset, called their {\em convex core}, see \cite{BS}. Denote by $\mathcal{GH}_{\Sigma,\theta}$ the space of convex GHM AdS metrics on $\Sigma\times\mathbb{R}$ with cone singularities of angles $\theta_{i}$ along the lines $\{p_{i}\}\times\mathbb{R}$, considered up to isotopies fixing each singular line (see the definition in Section 2.2).

    \subsection{Foliations of AdS spacetimes by $K$-surfaces}

        Our main result (Theorem \ref{foliation of the GHMC AdS manifold} below) asserts that in any convex GHM AdS spacetime with particles, the complement of the convex core admits a unique foliation by constant Gauss curvature surfaces. This extends to spacetimes with particles a result of B\'eguin, Barbot and Zeghib \cite{BBZ2} for non-singular GHM AdS spacetimes.

\begin{theorem}\label{foliation of the GHMC AdS manifold}
            Let $(N,g)$ be a convex GHM AdS manifold with particles and let $C(N)$ be the convex core of $N$. Then $N\setminus C(N)$ admits a unique foliation by locally strictly convex, constant Gauss curvature surfaces which are orthogonal to the singular lines.
\end{theorem}

    \subsection{Parameterization of the space of GHM AdS spacetimes}
        It is known that $\mathcal{GH}_{\Sigma,\theta}$ can be parameterized in several ways, such as the extension of Mess parametrization by $\mathcal{T}_{\Sigma,\theta}\times \mathcal{T}_{\Sigma,\theta}$ in terms of the left and right metrics, and the parametrization by $\mathcal{T}_{\Sigma,\theta}\times\mathcal{ML}_{\Sigma,n_0}$ in terms of the embedding data (the induced metric and the bending lamination) of the past (or future) boundary of the convex core. The first parametrization is equivalent to Thurston's Earthquake Theorem for hyperbolic metrics on $\Sigma$ with cone singularities of fixed angles less than $\pi$ (see \cite[Theorem 1.2]{BS}).

        Moreover, $\mathcal{GH}_{\Sigma,\theta}$ can also be parameterized by the cotangent bundle $T^{*}\mathcal{T}_{\Sigma,\theta}$ of $\mathcal{T}_{\Sigma,\theta}$, since 
        $T^{*}\mathcal{T}_{\Sigma,\theta}$ is homeomorphic to the quotient of the space $\mathcal{H}_{\Sigma,\theta}$ of maximal surfaces in germs of AdS manifolds with particles by diffeomorphisms isotopic to the identity fixing each marked point of $\Sigma$
        (see \cite[Theorem 5.11]{KS}) and there is a bijection between 
        this quotient space and $\mathcal{GH}_{\Sigma,\theta}$ (see \cite[Theorem 1.4]{Toulisse1}).

        We give a new parametrization of $\mathcal{GH}_{\Sigma,\theta}$ by $\mathcal{T}_{\Sigma,\theta}\times \mathcal{T}_{\Sigma,\theta}$ in terms of constant Gauss curvature surfaces. Specifically, we consider the map $\phi_{K}:\mathcal{T}_{\Sigma,\theta}\times \mathcal{T}_{\Sigma,\theta}\rightarrow \mathcal{GH}_{\Sigma,\theta}$, for each $K<-1$, which assigns to an element $(\tau,\tau')\in\mathcal{T}_{\Sigma,\theta}\times \mathcal{T}_{\Sigma,\theta}$ the isotopy class of the (unique) convex GHM AdS manifold $(N,g)$ with particles, such that it contains a future-convex, spacelike, constant curvature $K$ surface which is orthogonal to the singular lines, with induced metric $I\in \tau$ and third fundamental form $III\in \tau'$.

        \begin{theorem}\label{parametrization map}
            For any $K\in(-\infty, -1)$ and $\theta=(\theta_1,\cdots, \theta_{n_0})\in (0,\pi)^{n_0}$, the map $\phi_{K}:\mathcal{T}_{\Sigma,\theta}\times \mathcal{T}_{\Sigma,\theta}\rightarrow \mathcal{GH}_{\Sigma,\theta}$ is a homeomorphism.
        \end{theorem}

         Furthermore, we find that this result provides a convenient tool to prove the existence and uniqueness of the foliation of the complement of the convex core in a convex GHM AdS manifold with particles by locally strictly convex constant (Gauss) curvature surfaces which are orthogonal to the singular lines.

         In the case of a non-singular 3-dimensional GHM Lorentzian manifold of constant curvature, the corresponding result about the foliation by constant Gauss curvature surfaces has been proved by Barbot, B\'{e}guin and Zeghib (see Theorem 2.1 in \cite{BBZ2}). For the existence part, the argument in \cite{BBZ2} depends on the construction of barriers (see Definition 3.1 in \cite{BBZ2}) and a barriers theorem of Gerhardt (see \cite{Gerhardt}) to find the surface of a given constant curvature from the barriers. Here by contrast, Theorem \ref{foliation of the GHMC AdS manifold}  is obtained as a consequence of Theorem \ref{parametrization map}, and we obtain a simpler approach to prove the existence of the foliation without using the barriers argument.

        \begin{remark}
            For convenience, constant Gauss curvature surfaces are called constant curvature surfaces, or simply $K$-surfaces, henceforth.
        \end{remark}

    \subsection{Landslides on hyperbolic surfaces with cone singularities}
        Finally, we use the results obtained on $K$-surfaces in GHM AdS spacetimes with particles to extend some recent results on the landslide flow (see \cite{BMS1,BMS2}) to hyperbolic surfaces with cone singularities of fixed angles less than $\pi$.

        Landslides are transformations of hyperbolic structures on a closed surface $\Sigma_g$ of genus $g\geq2$, introduced in \cite{BMS1,BMS2} as ``smoother'' analogs of earthquakes. Earthquakes depend on the choice of a measured lamination $\lambda\in \cML_{\Sigma_g}$, so the earthquake flow can be defined as a map $$ \cE: \cT_{\Sigma_g}\times \cML_{\Sigma_g}\times \R\to \cT_{\Sigma_g}~ $$
        $$(h,\lambda,t)\mapsto E_{t\lambda}(h)~,$$
which for fixed $\lambda\in \cML_{\Sigma_g}$ defines an action of $\R$ on $\cT_{\Sigma_g}$.

        Landslide transformations, on the other hand, can be described as an action of $S^1$ on $\cT_{\Sigma_g}\times \cT_{\Sigma_g}$. For $e^{i\alpha}\in S^1$ and $(h,h')\in \cT_{\Sigma_g}\times \cT_{\Sigma_g}$, we denote by $L_{e^{i\alpha}}(h,h')\in \cT_{\Sigma_g}\times \cT_{\Sigma_g}$ the image of $(h,h')$ by the landslide flow, and $L^1_{e^{i\alpha}}(h,h')$ its projection on the first factor. If $(t_n)_{n\in \N}$ and $(h'_n)_{n\in \N}$ are sequences in $\R_{>0}$ and $\cT_{\Sigma_g}$, respectively, such that $t_nh'_n\to \lambda\in \cML_{\Sigma_g}$, then $L^1_{e^{it_n}}(h,h'_n)\to E_{\lambda/2}(h)$ as $n\to \infty$, see \cite[Theorem 1.12]{BMS1}.

        In Section 5, we use Theorem \ref{foliation of the GHMC AdS manifold} and other tools to extend the definition of landslide transformations to hyperbolic surfaces with cone singularities of fixed angles less than $\pi$. We show that the analog of Thurston's Earthquake Theorem extends to landslides on those hyperbolic cone surfaces: for all $h_1, h_2\in \cT_{\Sigma, \theta}$ and all $e^{i\alpha}\in S^1\setminus \{ 1\}$, there exists a unique $h'_1\in \cT_{\Sigma, \theta}$ such that $L^1_{e^{i\alpha}}(h_1,h'_1)=h_2$, see Theorem \ref{tm:landslides}.

        We then go on to deduce from the properties of the landslide flow further results on  the induced metrics and third fundamental forms of different $K$-surfaces in a given GHM AdS spacetime with particles, see Theorem \ref{tm:K1K2}.

    \subsection{Outline of the paper.}

        Section 2 is devoted to background material on different notions necessary for the rest of the paper: hyperbolic surfaces with cone singularities, AdS spacetimes with particles, etc. In Section 3 we prove Theorem \ref{parametrization map}, while Section 4 is devoted to the proof of Theorem \ref{foliation of the GHMC AdS manifold}. Section 5 describes applications to the landslide flow on the space of hyperbolic metrics with cone singuarities (of fixed angles) on a surface.


    \section{Background material}

    \subsection{Hyperbolic metrics with cone singularities.}
    First we recall the local model of a hyperbolic metric with a cone singularity of angle $\theta_{0}$.

    Let $\mathbb{H}^{2}$ be the Poincar\'{e} model of the hyperbolic plane. Denote by $\mathbb{H}^{2}_{\theta_{0}}$ the space obtained by taking a wedge of angle $\theta_{0}$ bounded by two half-lines intersecting at the center 0 of $\mathbb{H}^{2}$ and gluing the two half-lines by a rotation fixing 0. We call $\mathbb{H}^{2}_{\theta_{0}}$ the \emph{hyperbolic disk with cone singularity of angle $\theta_{0}$}, which is a punctured disk with the induced metric
    \begin{equation*}
        g_{\theta_{0}}=dr^{2}+\sinh^{2}(r)d\alpha^{2},
    \end{equation*}
    where $(r,\alpha)\in\mathbb{R}_{>0}\times\mathbb{R}/\theta_{0}\mathbb{Z}$ is a polar coordinate of $\mathbb{H}^{2}_{\theta_{0}}$.

    In conformal coordinates, $g_{\theta_{0}}$ has the following expression, which is obtained by pulling back the Poincar\'{e} metric by the map $z\mapsto z^{t}/t$ with $t=\theta_{0}/2\pi$.
    \begin{equation*}
        g_{\theta_{0}}=\frac{4|z|^{2(t-1)}}{(1-t^{-2}|z|^{2t})^2}|dz|^{2}.
    \end{equation*}

    To apply the existence of harmonic maps between Riemann surfaces with marked points and hyperbolic surfaces with cone singularities (see \cite[Theorem 2]{GR}) in subsequent sections , we need a regularity condition of the metric around the cone singularities and we introduce the following weighted H\"older spaces (see \cite[Section 2.2]{GR} and \cite[Definition 2.1]{Toulisse2}).

    \begin{definition}
        For $R>0$, let $D(R):=\{z\in\mathbb{C},|z|\in(0,R)\}$. A function $f:D(R)\rightarrow\mathbb{C}$ is said to be in $\chi^{0,\gamma}_b(D(R))$ with $\gamma\in(0,1)$ if
        \begin{equation*}
            ||f||_{\chi^{0,\gamma}_{b}}:= \sup_{z\in D(R)}|f(z)|+\sup_{z,z'\in D(R)} \frac{|f(z)-f(z')|}
            {|\alpha-\alpha'|^{\gamma}+|\frac{r-r'}{r+r'}|^{\gamma}}
            <\infty,
        \end{equation*}
        where $z=re^{i\alpha}$ and $z'=r'e^{i\alpha'}$. Let $k\in\mathbb{N}$, we say that $f\in\chi^{k,\gamma}_b (D(R))$ if $(r\partial_r)^i\partial^{j}_{\alpha}f$ is in ${\chi^{0,\gamma}_{b}}(D(R))$ for all $i+j\leq k$. In particular, this implies that $f\in \mathcal{C}^{k}(D(R))$.
    \end{definition}

    \begin{definition}
        Let $\mathfrak{p}=\{p_1,...,p_{n_0}\}$, $\theta=\{\theta_1,...,\theta_{n_{0}}\}$ and $\Sigma_{\mathfrak{p}}=\Sigma\setminus \mathfrak{p}$. A hyperbolic metric on $\Sigma$ with cone singularities of angle $\theta$ at $\mathfrak{p}$ is a (singular) metric $g$ on $\Sigma$ with the property: for each compact set $K\subset \Sigma_{\mathfrak{p}}$, $g|_K$ is $\mathcal{C}^{2}$ and has constant curvature $-1$, and for each marked point $p_{i}$,  there exist a neighbourhood $U_{i}$ with local conformal coordinates $z$ centered at $p_i$ and a local diffeomorphism $\psi\in\chi^{2,\gamma}_{b}(U_i)$ such that $g|_{U_i}$ is the pull back by $\psi$ of the metric $g_{\theta_i}$. Denote by $\mathfrak{M}^{\theta}_{-1}$ the space of hyperbolic metrics on $\Sigma$ with cone singularities of angle $\theta$ at $\mathfrak{p}$.
    \end{definition}

    We say that $f$ is a \emph{diffeomorphism} of $\Sigma_{\mathfrak{p}}$ if for each compact set $K\subset \Sigma_{\mathfrak{p}}$, $f|_{K}$ is of class $\mathcal{C}^{3}$ and for each marked point $p_i$, there exists a neighbourhood $U_i$ of $p_i$ such that $f|_{U_i}\in\chi^{2,\gamma}_{b}(U_i)$. Denote by $\mathfrak{Diff}_{0}(\Sigma_\mathfrak{p})$ the space of diffeomorphisms on $\Sigma_\mathfrak{p}$ which are isotopic to the identity (fixing each marked point). They act by pull-back on $\mathfrak{M}^{\theta}_{-1}$. We say that two metrics $h_1,h_2\in\mathfrak{M}^{\theta}_{-1}$ are \emph{isotopic} if there exists a map $f\in\mathfrak{Diff}_{0}(\Sigma_\mathfrak{p})$ such that $h_1$ is the pull back by $f$ of $h_2$.

    Denote by $\mathcal{T}_{\Sigma,\theta}$ the \emph{Teichm\"uller space of hyperbolic metrics on $\Sigma$ with fixed cone angle $\theta$}, which is the space of isotopy classes of hyperbolic metrics on $\Sigma$ with cone singularities of angle $\theta$ at $\mathfrak{p}$. Note that $\mathcal{T}_{\Sigma,\theta}=\mathfrak{M}^{\theta}_{-1}/\mathfrak{Diff}_{0}(\Sigma_\mathfrak{p})$ and $\mathfrak{M}^{\theta}_{-1}$ is a differentiable submanifold of the manifold consisting of all $\mathcal{H}^2$ symmetric (0,2)-type tensor fields.
    $\mathcal{T}_{\Sigma,\theta}$ is a finite-dimensional differentiable manifold which inherits a natural quotient topology.

    \subsection{Convex GHM AdS manifolds with particles.}

    First we recall the related notations and terminology in order to define convex GHM AdS manifolds with particles.

    \textbf{The AdS 3-space.} Let $\mathbb{R}^{2,2}$ be $\mathbb{R}^4$ with the quadratic form $q(x)=x_{1}^{2}+x_{2}^{2}-x_{3}^{2}-x_{4}^{2}$. The \emph{anti-de Sitter (AdS) 3-sapce} is defined as the quadric:
    \begin{equation*}
        {AdS}_{3}=\{x\in\mathbb{R}^{2,2}:q(x)=-1\}.
    \end{equation*}
    It is a 3-dimensional Lorentzian symmetric space of constant curvature $-1$ diffeomorphic to $\mathbb{D}\times S^{1}$, where $\mathbb{D}$ is a 2-dimensional disk.

    Consider the projective map $\pi:\mathbb{R}^{2,2}\backslash\{0\}\rightarrow \mathbb{RP}^{3}$. The \emph{Klein model} $\mathbb{ADS}_{3}$ of AdS 3-space is defined as the image of ${AdS}_{3}$ under the projection $\pi$. It is clear that $\mathbb{ADS}_{3}=\pi({AdS}_{3})={AdS}_{3}/\{\pm id\}$. The boundary $\partial \mathbb{ADS}_{3}$ is the image of the quadratic $Q=\{x\in\mathbb{R}^{2,2}:q(x)=0\}$ under $\pi$, which is foliated by two families of projective lines, called the left and right leaves, respectively.

    Geodesics in the Klein model $\mathbb{ADS}_{3}$ are given by projective lines: the spacelike geodesics correspond to the projective lines intersecting the boundary $\partial \mathbb{ADS}_{3}$ in two points, while lightlike geodesics are tangent to $\partial \mathbb{ADS}_{3}$, and timelike geodesics do not intersect $\partial \mathbb{ADS}_{3}$.

    The group Isom$_{0}(\mathbb{ADS}_{3})$ of space and time orientation preserving isometries of $\mathbb{ADS}_{3}$ can be identified as PSL(2, $\mathbb{R})\times$ PSL(2, $\mathbb{R}$).

    \textbf{The singular AdS 3-space.} Let $\theta_0>0$. Define the \emph{singular AdS 3-space} of angle $\theta_0$ as
    \begin{equation*}
        AdS^3_{\theta_0}:=\{(t,r, \alpha)\in\mathbb{R}\times \mathbb{R}^{+}\times\mathbb{R}/\theta_{0}\mathbb{Z}\}.
    \end{equation*}
    with the metric
    \begin{equation*}
        -dt^{2}+\cos^{2}t(dr^{2}+\sinh^{2}(r)d\alpha^{2}).
    \end{equation*}
    The set corresponding to $r=0$ is called the \emph{singular line} in $AdS^3_{\theta_0}$.

    It is clear that $AdS^3_{\theta_0}$ is a Lorentzian manifold of constant curvature $-1$ outside the singular line, that is, it is locally modelled on the universal cover of $AdS_3$. Indeed, it is obtained from the complete hyperbolic surface with a cone singularity of angle $\theta_{0}$ by taking a warped product with $\mathbb{R}$ (see e.g.\cite{BS,KS,Toulisse1}).

    An embedded surface in $AdS^3_{\theta_0}$ is \emph{spacelike} if it intersects the singular line at exactly one point and it is spacelike outside the intersection with the singular locus.

    \textbf{AdS manifolds with particles.} An \emph{AdS manifold with particles} is a (singular) Lorentzian 3-manifold in which any point $x$ has a neighbourhood isometric to a subset of $AdS^3_{\theta_0}$ for some $\theta_{0}\in(0,\pi)$.

    A closed embedded surface $S$ in an AdS manifold with particles is \emph{spacelike} if it is locally modelled on a spacelike surface in $AdS^3_{\theta_0}$ for some $\theta_0\in(0,\pi)$.

    \begin{definition}
        Let $S\subset AdS^3_{\theta_0}$ be a spacelike surface which intersects the singular line at a point $x$. $S$ is orthogonal to the singular locus at $x$ if the distance to the totally geodesic plane $P$ orthogonal to the singular line at $x$ satisfies:
        \begin{equation*}
            \lim_{y\rightarrow x,y\in S}\frac{d(y,P)}{d_S(x,y)}=0,
        \end{equation*}
        where $d_S(x,y)$ is the distance between $x$ and $y$ along $S$.

        If now $S$ is a spacelike surface in an AdS manifold $M$ with particles which intersects a singular line $l$ at a point $x'$. $S$ is said to be orthogonal to $l$ at $x'$ if there exists a neighborhood $U\subset M$ of $x'$ which is isometric to a neighborhood of a singular point in $AdS^3_{\theta_0}$ such that the isometry sends $S\cap U$ to a surface orthogonal to the singular line in $AdS^3_{\theta_0}$.
    \end{definition}

    \begin{definition}
        An AdS manifold $M$ with particles is convex GHM if
        \begin{itemize}
            \item $M$ is convex GH: it contains a locally convex spacelike surface $S$ orthogonal to the singular lines, which intersects every inextensible timelike curve exactly once.
            \item $M$ is maximal: if any isometric embedding of $M$ into a convex GH AdS manifold is an isometry.
        \end{itemize}
    \end{definition}

    Let $\mathcal{GH}'_{\Sigma,\theta}$ be the space of convex GHM AdS metrics on $\Sigma\times \mathbb{R}$ with cone singularities of angles $\theta_i$ along the lines $\{p_i\}\times \mathbb{R}$. Denote by $\mathfrak{Diff}_0(\Sigma\times\mathbb{R})$ the space of diffeomorphisms on $\Sigma\times\mathbb{R}$ isotopic the identity fixing each singular line. We say that two metrics $g_1,g_2\in\mathcal{GH}'_{\Sigma,\theta}$ are \emph{isotopic} if there exists a map $f\in\mathfrak{Diff}_0(\Sigma\times\mathbb{R})$ such that $g_1$ is the pull back by $f$ of $g_2$.

    Denote by $\mathcal{GH}_{\Sigma,\theta}$ the \emph{space of convex GHM AdS metrics on $\Sigma\times \mathbb{R}$ with particles of fixed angle $\theta$},which is the space of isotopy classes of convex GHM AdS metrics with cone singularities of angles $\theta_i$ along the lines $\{p_i\}\times \mathbb{R}$. Note that $\mathcal{GH}_{\Sigma,\theta}=\mathcal{GH}'_{\Sigma,\theta}/\mathfrak{Diff}_0(\Sigma\times\mathbb{R})$ and it is a finite-dimensional differentiable manifold with a natural quotient topology.

    \subsection{Convex spacelike surfaces in a convex GHM AdS manifold with particles.}

    Let $(N,g)$ be a convex GHM AdS manifold with particles.
    Let $S\subset N$ be an (embedded) spacelike surface orthogonal to the singular lines with the induced metric $I$. The shape operator $B:TS\rightarrow TS$ of $S$ is defined as
    \begin{equation*}
        B(u)=\nabla_{u}n,
    \end{equation*}
    where $n$ is the future-directed unit normal vector field on $S$ and $\nabla$ is the Levi-Civita connection of $(N,g)$. The second and third fundamental forms of $S$ are defined respectively as
    \begin{equation*}
        II(u,v)=I(Bu,v),\qquad
        III(u,v)=I(Bu,Bv).
    \end{equation*}

    \begin{definition}
        Let $S$ be a convex spacelike surface orthogonal to the singular lines in a convex GHM AdS manifold $N$ with particles. We say that $S$ is future-convex (resp. past-convex) if its future $I^{+}(S)$ (resp. its past $I^{-}(S)$) is geodesically convex. We say that $S$ is strictly future-convex (resp. strictly past-convex) if $I^{+}(S)$ (resp. $I^{-}(S)$) is strictly geodesically convex.
    \end{definition}

    Note that if $S$ is future-convex (resp. past-convex), then for each regular point $x$ of $S$, both the principal curvatures at $x$ are non-negative (resp. non-positive). If $S$ is strictly future-convex (resp. strictly past-convex), then for each regular point $x$ of $S$, both the principal curvatures at $x$ are positive (resp. negative).

    \subsection{The duality between strictly convex surfaces in convex GHM AdS manifolds with particles}

    First we recall the duality between points and hyperplanes in $AdS_{3}$ (see e.g. \cite{BMS2,BBZ2}).

    Observe that $AdS_3$ is a quadric in $\mathbb{R}^{2,2}$. Every point $x$ in
    $AdS_3$ is exactly the intersection in $\mathbb{R}^{2,2}$ of $AdS_3$ with a ray $l$ starting from the origin 0 on which the quadratic form is negative definite.
    Denote by $l^\perp$ the hyperplane orthogonal to $l$ in $\mathbb{R}^{2,2}$, with the induced metric of signature (2,1). The intersection between $l^\perp$ and $AdS_3$ is the disjoint union of two totally geodesic spacelike planes $P^{\pm}_{x}$, where $P^{+}_{x}$ (resp. $P^{-}_{x}$) is at a distance $\pi/2$ in the future (resp. in the past) of $x$.

    Conversely, every totally geodesic spacelike plane $P$ in $AdS_3$ is the intersection of $AdS_3$ with a hyperplane $H$ of signature (2,1) in $\mathbb{R}^{2,2}$. The orthogonal $H^{\perp}$ of $H$ intersects $AdS_3$ at two antipodal points $x^{\pm}_{P}$, where $x^{+}_{P}$ (resp. $x^{-}_{P}$) is at a distance $\pi/2$ in the future (resp. in the past) of $P$.

    We define the dual $P^{*}$ of $P$ as the past (resp. future) intersection point $x^{-}_{P}$ if the dual $x^{*}$ of $x$ is defined to be $AdS_3\cap P^{+}_{x}$ (resp. $AdS_3\cap P^{-}_{x}$). The \emph{dual surface} $S^{*}$ of a strictly convex surface $S\subset AdS_{3}$ is defined as the set of points on the convex side of $S$ which are the dual points of the support planes of $S$. Equivalently, $S^{*}$ can be obtained by pushing $S$ along orthogonal geodesics on the convex side for a distance $\pi/2$ (see \cite[Proposition 11.9]{BBZ2}).

    Note that $AdS^3_{\theta_0}$ can be obtained from the universal cover of $AdS_3$ by taking a wedge of angle $\theta_0$ bounded by two timelike totally geodesic half-planes and gluing the two half-planes by a rotation fixing the common timelike geodesic. For a strictly convex spacelike surface $S\subset AdS^3_{\theta_0}$ orthogonal to the singular lines, there is a natural generalization for the dual surface $S^*$.

    Since an AdS manifold with particles is locally modelled on $AdS^3_{\theta_0}$ for some $\theta_0\in(0,\pi)$, we can generalize to the singular case the duality between strictly convex spacelike surfaces.

    \begin{definition}
        Let $S$ be a strictly convex spacelike surface orthogonal to the singular lines in a convex GHM AdS manifold $N$ with particles. The \emph{dual surface} $S^{*}$ of $S$ is defined as the surface obtained by pushing $S$ along  orthogonal geodesics on the convex side for a distance $\pi/2$.
    \end{definition}

    Observe that the surface obtained by pushing a strictly convex surface $S\subset N$ (which is orthogonal to the singular lines) along orthogonal geodesics on the convex side for a distance $t\in [0,\pi/2]$ is still orthogonal to the singular lines. The relation between dual strictly convex surfaces in GHMC AdS manifolds (see e.g.\cite{BBZ2,BMS2}) can be directly generalized to the following case with cone singularities.

    \begin{proposition}
        \label{the duality between convex surfaces}
        Let $(N,g)$ be a convex GHM AdS manifold with particles. Assume that $S\subset N$ is a strictly convex spacelike surface of constant curvature $K$ orthogonal to singular lines. Then
        \begin{enumerate}[(1)]
            \item $S^{*}$ is a strictly convex spacelike surface of constant curvature $K^{*}$ with the shape operator of opposite definiteness, which is orthogonal to the singular lines in $N$, where $K^{*}=-K/(1+K)$.
            \item The pull back of the induced metric on $S^{*}$ through the duality map is the third fundamental form of $S$ and vice versa.
            \item The dual surface $(S^{*})^{*}$ of $S^{*}$ is exactly $S$.
        \end{enumerate}
    \end{proposition}

    \subsection{Minimal Lagrangian maps between hyperbolic surfaces with cone singularities.}

    The construction of the parametrization of $\mathcal{GH}_{\Sigma,\theta}$ here depends strongly on minimal Lagrangian maps between hyperbolic surfaces with cone singularities.

    \begin{definition}
        Given two hyperbolic metrics $h,h'$ on $\Sigma$ with cone singularities. A \emph{minimal Lagrangian map} $m:(\Sigma,h)\rightarrow(\Sigma,h')$ is an area-preserving and orientation-preserving diffeomorphism such that its graph is a minimal surface in $(\Sigma\times\Sigma,h\oplus h')$.
    \end{definition}

    The following result is \cite[Theorem 1.3]{Toulisse1}.

    \begin{theorem}[Toulisse]
        \label{Toulisse1}
        Let $h,h'\in \mathfrak{M}^{\theta}_{-1}$. Then there exists a unique minimal Lagrangian diffeomorphism $m:(\Sigma,h)\rightarrow (\Sigma,h')$ isotopic to the identity.
    \end{theorem}

    This is shown by proving the existence and uniqueness of maximal surfaces (see \cite[Theorem 1.4]{Toulisse1}) in a convex GHM AdS manifold with particles.

    Recall that a spacelike surface of a convex GHM AdS manifold $(N,g)$ with particles  is said to be a \emph{maximal surface} if it is a locally area-maximizing Cauchy surface which is orthogonal to the singular lines.
    In particular, it has everywhere vanishing mean curvature and its principal curvatures are everywhere in $(-1,1)$ (see \cite[Lemma 5.15]{KS})
     and tend to zero at the intersections with particles (see \cite[Proposition 3.7]{Toulisse1}). It is described in \cite[Definition 5.10]{KS}
     that the space $\mathcal{H}_{\Sigma,\theta}$ of maximal surfaces in germs of AdS manifolds with particles has the following convenient properties.

    \begin{lemma}
        \label{maximal surface in germs}
        The space $\mathcal{H}_{\Sigma,\theta}$ is identified with the space of couples $(g,h)$, where $g$ is a smooth metric on $\Sigma_\mathfrak{p}$ with cone singularities of angle $\theta_i$ at the marked points $p_{i}$ for $i=1,...,n_0$ and $h$ is a symmetric bilinear form on $T\Sigma$ defined outside the marked points, such that
        \begin{itemize}
            \item $tr_{g}(h)=0$.
            \item $d^{\nabla}B=0$, where $\nabla$ is the Levi-Civita connection of $g$.
            \item $K_{g}=-1-\det_{g}(h)$.
        \end{itemize}
    \end{lemma}

    For the convenience of computation, we also introduce the following proposition, see \cite[Proposition 3.12]{KS}.
    \begin{proposition}
        \label{computation of connection and curvature}
        Let $\Sigma$ be a surface with a Riemann metric $g$. Let $A:T\Sigma\rightarrow T\Sigma$ be a bundle morphism such that $A$ is everywhere invertible and $d^{\nabla}A=0$, where $\nabla$ is the Levi-Civita connection of $g$. Let $h$ be the symmetric (0,2)-tensor defined by $h=g(A\bullet,A\bullet)$. Then the Levi-Civita connection of $h$ is given by
        \begin{equation*}
            \nabla^{h}_{u}(v)=A^{-1}\nabla_{u}(Av),
        \end{equation*}
         and its curvature is given by
         \begin{equation*}
            K_{h}=\frac{K_{g}}{\det(A)}.
         \end{equation*}
    \end{proposition}

    Minimal Lagrangian maps between hyperbolic surfaces with metrics in $\mathfrak{M}^{\theta}_{-1}$ have an equivalent description in terms of morphisms between tangent bundles (see e.g.\cite[Proposition 6.3]{Toulisse1}).

    \begin{proposition}
        \label{Toulisse2}
        Let $h,h'\in \mathfrak{M}^{\theta}_{-1}$. Then $m:(\Sigma,h)\rightarrow (\Sigma,h')$ is a minimal Lagrangian map if and only if there exists a bundle morphism $b:T\Sigma\rightarrow T\Sigma$ defined outside the singular locus which satisfies the following properties:
        \begin{itemize}
            \item $b$ is self-adjoint for $h$ with positive eigenvalues.
            \item $\det(b)=1$.
            \item $b$ satisfies the Codazzi equation: $d^{\nabla}b=0$, where $\nabla$ is the Levi-Civita connection of $h$.
            \item $h(b\bullet,b\bullet)$ is the pull back of $h'$ by a diffeomorphism $m:\Sigma \rightarrow\Sigma$ fixing each marked point.
            \item Both eigenvalues of $b$ tend to 1 at the cone singularities.
        \end{itemize}
    \end{proposition}

    \begin{proof}
        Note that Proposition 6.3 in \cite{Toulisse1} provides the equivalence between the existence of a minimal Lagrangian map $m:(\Sigma,h)\rightarrow (\Sigma,h')$ and the existence of a bundle morphism $b$ satisfies the first three properties. It suffices to check that given a minimal Lagrangian map $m:(\Sigma,h)\rightarrow (\Sigma,h')$, the bundle morphism $b$ also satisfies the last property. Set
        \begin{equation*}
            I'=\frac{1}{4}h((E+b)\bullet,(E+b)\bullet).
        \end{equation*}
        Denote by $J'$ the complex structure of $I'$ and set
        \begin{equation*}
         B'=-J'(E+b)^{-1}(E-b).
        \end{equation*}
        Moreover, $J'=(E+b)^{-1}J(E+b)$, where $J$ is the complex structure of $h$.

        Note that $J'B'=(E+b)^{-1}(E-b)$. It is not hard to check that $B'$ satisfies the following conditions:
        \begin{itemize}
            \item $B'$ is self-adjoint for $I'$. Indeed, choosing a suitable basis such that 
                $b$ 
                is diagonal and using the fact that $\det(b)=1$, we have
                \begin{equation*}
                    tr(J'B')=0,
                \end{equation*}
                which implies that $B'$ is self-adjoint for $I'$.
            \item $tr(B')=0$. This follows from the fact that $J'B'$ is self-adjoint for $I'$, since $E\pm b$ is self-adjoint for $h$ and $E+b$ commutes with $E-b$.
            \item $d^{\nabla^{'}}B'=0$, where $\nabla^{'}$ is the Levi-Civita connection of $I'$. Indeed, by Theorem \ref{computation of connection and curvature} and direct computation, we get
                \begin{equation*}
                    d^{\nabla^{'}}(J'B')=(E+b)^{-1}d^{\nabla}(E-b)=0.
                \end{equation*}
                Note that $J'$ is parallel for $\nabla^{'}$, it follows that $d^{\nabla^{'}}B'=0$.
            \item $K_{I'}=-1-\det(B')$. Indeed, by computation, we have
                \begin{equation*}
                    E+J'B'=2(E+b)^{-1}.
                \end{equation*}
            By Proposition \ref{computation of connection and curvature}, it follows that
            \begin{equation*}
                K_{I'}
                =\frac{K_{h}}{{\det}(\frac{1}{2}(E+b))}
                =-{\det}(2(E+b)^{-1})
                =-{\det}(E+J'B')=-1-{\det}(B').
            \end{equation*}
        \end{itemize}

        Set $II'=I'(B'\bullet,\bullet)$. By Lemma \ref{maximal surface in germs}, the couple $(I',II')$ is exactly the first and second fundamental form of a maximal surface $S'$ in a convex GHM AdS manifold $(N',g')$ with particles, where $(N',g')$ has the left metric
        \begin{equation*}
            I'((E+J'B')\bullet,(E+J'B')\bullet)=h,
        \end{equation*}
        and the right metric
        \begin{equation*}
            I'((E-J'B')\bullet,(E-J'B')\bullet)=h'.
        \end{equation*}

        Note that the eigenvalues of $B'$ tend to zero at the intersections of $S'$ with the particles (see \cite[Proposition 3.7]{Toulisse1}) and $B'=-J'(E+b)^{-1}(E-b)$. Then both eigenvalues of $b$ tend to 1 at the cone singularities. This completes the proof.
    \end{proof}

    \begin{corollary}\label{hyperbolic metrics and bundle morphism}
        Let $h,h'\in \mathfrak{M}^{\theta}_{-1}$. Then there exists a unique bundle morphism $b:T\Sigma\rightarrow T\Sigma$ defined outside the singular locus, which is self-adjoint for $h$ with positive eigenvalues, has determinant 1 and satisfies the Codazzi equation: $d^{\nabla}b=0$, where $\nabla$ is the Levi-Civita connection of $h$, such that $h(b\bullet,b\bullet)$ is isotopic to $h'$. Moreover, both eigenvalues of $b$ tend to 1 at the cone singularities.
    \end{corollary}

    \begin{definition}\label{defintion of normalized representatives}
        We say that a pair of hyperbolic metrics $(h,h')$ is normalized if there exists a bundle morphism $b:T\Sigma\rightarrow T\Sigma$ defined outside the singular locus, which is self-adjoint for $h$, has determinant $1$, and satisfies the Codazzi equation, such that $h'=h(b\bullet,b\bullet)$, or equivalently if the identity from $(\Sigma,h)$ to $(\Sigma,h')$ is a minimal Lagrangian diffeomorphism.
    \end{definition}

    \begin{remark}\label{rk: normalized representatives}
        By Corollary \ref{hyperbolic metrics and bundle morphism}, for any $(\tau,\tau')\in\mathcal{T}_{\Sigma,\theta}\times \mathcal{T}_{\Sigma,\theta}$, we can realize $(\tau,\tau')$ as a normalized representative $(h,h')$. Note that the normalized representative of $(\tau,\tau')$ is unique up to isotopies acting diagonally on both $h$ and $h'$.
    \end{remark}

    \section{Parametrization of $\mathcal{GH}_{\Sigma,\theta}$ in terms of constant curvature surfaces.}

    \subsection{The definition of the map $\phi_{K}$.}
    For the construction of the map $\phi_{K}$, we introduce the following proposition which ensures the existence and the uniqueness (up to isometries) of the maximal extension of a convex GH AdS manifold with particles (see \cite[Proposition 2.6]{BS}).

    \begin{proposition}
        \label{uniqueness of maximal extension}
            Let $(N,g)$ be a convex GH AdS manifold with particles. There exists a unique (considered up to isometries) convex GHM AdS manifold $(N',g')$ with particles, called the maximal extension of $(N,g)$, in which $(N,g)$ can be isometrically embedded.
    \end{proposition}

    \begin{lemma}
        \label{definition of the map}
            Let $K\in(-\infty, -1)$ and let $(h,h')\in\mathfrak{M}^{\theta}_{-1}\times\mathfrak{M}^{\theta}_{-1}$ be a pair of normalized metrics. Then there exists a unique GHM AdS manifold $(N,g)$ that contains a future-convex, spacelike, constant curvature $K$ surface which is orthogonal to the singular lines, with the induced metric $I=(1/|K|)h$ and the third fundamental form $III=(1/|K^{*}|)h'$, where $K^{*}=-K/(1+K)$.
    \end{lemma}

    \begin{proof}
        Let $b:T\Sigma\rightarrow T\Sigma$ be the bundle morphism associated to $h$ and $h'$ by Definition \ref{defintion of normalized representatives}, so that $h'=h(b\bullet, b\bullet)$.

        Let $I=(1/|K|)h$. We equip $\Sigma$ with the metric $I$ and consider a bundle morphism $B:T\Sigma\rightarrow T\Sigma$, which is defined by $B=\sqrt{-1-K}b$. By the properties of $h$ and $b$, it follows that
        \begin{itemize}
            \item $(\Sigma, I)$ has constant curvature $K$.
            \item $B$ is self-adjoint for $I$ with positive eigenvalues.
            \item $B$ satisfies the Codazzi equation: $d^{\nabla^{I}}B=0$, where $\nabla^{I}$ is the Levi-Civita connection of $I$.
            \item $B$ satisfies the Gauss equation: $K=-1-\det(B)$.
        \end{itemize}

        Consider the manifold $\Sigma\times [0,\frac{\pi}{2})$ with the following metric:
        \begin{equation*}
            g_{0}=-dt^{2}+I((\cos(t)E+\sin(t)B)\bullet,(\cos(t)E+\sin(t)B)\bullet)~,
        \end{equation*}
        where $E$ is the identity isomorphism on $T\Sigma$ and $t\in[0,\frac{\pi}{2})$. Note that for each $t\in [0,\frac{\pi}{2})$, the surface $\Sigma\times \{t\}$ is the equidistant surface at distance $t$ from the surface $\Sigma\times\{0\}$ on the convex side. The Lorentzian metric $g_{0}$ is a convex GH AdS metric on $\Sigma\times[0,\frac{\pi}{2})$ with cone singularities of angle $\theta_{i}$ along the line $\{p_{i}\}\times [0,\frac{\pi}{2})$.

        By Proposition \ref{uniqueness of maximal extension}, there exists a unique maximal extension $(N,g)$ of the AdS manifold $(\Sigma\times [0,\frac{\pi}{2}),g_{0})$ with particles, which is a convex GHM AdS manifold with particles, such that the restriction of $g$ to the subset $\Sigma\times [0,\frac{\pi}{2})$ of $N$ is exactly $g_{0}$.

        Since $B$ has positive eigenvalues, the embedded surface $\Sigma\times \{ 0\}$ is future-convex. Hence, $N$ contains a future-convex, spacelike, constant curvature $K$ surface which is orthogonal to the singular lines, with the induced metric $I=(1/|K|)h$ and the third fundamental form
        \begin{equation*}
            III=I(B\bullet,B\bullet)
            =\frac{1}{|K|}h(\sqrt{-1-K}b\bullet,\sqrt{-1-K}b\bullet)
            =\frac{1}{|K^{*}|}h',
        \end{equation*}
        where $|K^{*}|=-K/(1+K)$. This shows the existence of the required manifold $(N,g)$.

        Now we show the uniquess of $(N,g)$. Suppose that $(N_{1},g_{1})$ is another convex GHM AdS manifold with particles which contains such a required surface $S_{1}$. Then $S_{1}$ has the induced metric $I_{1}=(1/|K|)h=I$ with shape operator $B_{1}$ and third fundamental form
        \begin{equation*}
            III_1=I(B_{1}\bullet,B_{1}\bullet)
            =\frac{1}{|K^{*}|}h'=I(B\bullet,B\bullet)
            =III.
        \end{equation*}

        Since $S_{1}$ is future-convex, then $B_{1}$ is positive definite. Therefore, the shape operator $B_{1}$ of $S_{1}$ in $(N_{1},g_{1})$ is equal to $B$. Note that the embedding data $(\Sigma, I, B)$ is exactly $(\Sigma, I_{1}, B_{1})$, then $(N_{1},g_{1})=(N,g)$. This completes the proof.
    \end{proof}

    \begin{lemma}\label{well-define of the map}
        For any $(\tau, \tau')\in\mathcal{T}_{\Sigma,\theta}\times \mathcal{T}_{\Sigma,\theta}$, let $(h,h')$ and $(h_{1}, h_{1}')$ be two normalized representatives of $(\tau,\tau')$. Let $(N,g)$ and $(N_{1},g_{1})$ be the convex GHM AdS manifolds with particles associated to $(h,h')$ and $(h_{1},h_{1}')$, as described in Lemma \ref{definition of the map}. Then $(N,g)$ is isotopic to $(N_{1},g_{1})$.
    \end{lemma}

    \begin{proof}
        Note that $(h,h')$, $(h_{1},h_{1}')$ are normalized representatives of $(\tau,\tau')$. By Remark \ref{rk: normalized representatives}, there exists a diffeomorphism $\varphi$ from $\Sigma$ to $\Sigma$ which is isotopic to the identity (the isotopy fixes the marked points), such that $h_{1}=\varphi^{*}h$ and $h_{1}'=\varphi^{*}h'$.

        Let $(\Sigma,I,B,III)$ and $(\Sigma,I_{1},B_{1},III_{1})$ be the corresponding data of the surface contained in $(N,g)$ and $(N_{1},g_{1})$, as described in Lemma \ref{definition of the map}, respectively. Then $I=(1/|K|)h$, $III=(1/|K^{*}|)h'$ and  $I_{1}=(1/|K|)h_{1}$, $III_{1}=(1/|K^{*}|)h_{1}'$.
        It follows that
        \begin{equation}
            \label{equation1}
                I_{1}=\varphi^{*}(I),\qquad  III_{1}=\varphi^{*}(III).
        \end{equation}

        To see $(N_{1},g_{1})$ is isotopic to $(N,g)$, it suffices to prove that $II_{1}=\varphi^{*}(II)$, where $II$ is the second fundamental form of $(\Sigma,I)$ in $(N,g)$ and $II_{1}$ is the second fundamental form of $(\Sigma,I_{1})$ in $(N_{1},g_{1})$.

         By \eqref{equation1}, we have
         \begin{equation*}
            III_{1}
            =\varphi^{*}(III)
            =\varphi^{*}(I(B\bullet, B\bullet))
            =\varphi^{*}(I(B^{2}\bullet,\bullet))
            =I(B^{2}d\varphi\bullet, d\varphi\bullet),
         \end{equation*}
         where $d\varphi$ denotes the differential map (or the Jacobian matrix) of $\varphi$.

         Note that
         \begin{equation*}
            \uppercase\expandafter{\romannumeral3}_{1}
            =\uppercase\expandafter{\romannumeral1}_{1}(B_{1}\bullet, B_{1}\bullet)=(\varphi^{*}I)(B_{1}\bullet,B_{1}\bullet)
            =(\varphi^{*}I)(B_{1}^{2}\bullet,\bullet)
            =I(d\varphi B_{1}^{2}\bullet,d\varphi\bullet)
         \end{equation*}

         Hence, $B^{2}=(d\varphi)B_{1}^{2}(d\varphi)^{-1}$. Denote $A=(d\varphi)B_{1}(d\varphi)^{-1}$ and hence $B^{2}=A^{2}$. Since both $A$ and $B$ are self-adjoint with positive eigenvalues, an elementary argument shows that $A=B$, that is, $(d\varphi)B_{1}=B(d\varphi)$.
         Therefore,
         \begin{equation*}
            \varphi^{*}(II)
            =\varphi^{*}(I(B\bullet,\bullet))
            =I(Bd\varphi\bullet,d\varphi\bullet)
            =I(d\varphi B_{1}\bullet,d\varphi\bullet)
            =(\varphi^{*}I)(B_{1}\bullet,\bullet)
            =I_{1}(B_{1}\bullet,\bullet)
            =\uppercase\expandafter{\romannumeral2}_{1}.
         \end{equation*}
         This completes the proof of Lemma \ref{well-define of the map}.
    \end{proof}

    \begin{definition}
        For any $K\in(-\infty, -1)$, define the map $\phi_{K}:\mathcal{T}_{\Sigma,\theta}\times \mathcal{T}_{\Sigma,\theta}\rightarrow \mathcal{GH}_{\Sigma,\theta}$ by assigning to an element $(\tau,\tau')\in\mathcal{T}_{\Sigma,\theta}\times \mathcal{T}_{\Sigma,\theta}$ the isotopy class of the convex GHM AdS manifold $(N,g)$ with particles satisfying the prescribed property in Lemma \ref{definition of the map}. As a consequence of Lemma \ref{definition of the map} and Lemma \ref{well-define of the map}, this map is well-defined.
    \end{definition}

    \begin{remark}
        For convenience, for each pair $(\tau,\tau')\in\mathcal{T}_{\Sigma,\theta}\times \mathcal{T}_{\Sigma,\theta}$, we always represent it by a pair of normalized hyperbolic metrics $(h,h')$ and represent $\phi_{K}(\tau,\tau')$ by the convex GHM AdS manifold $(N,g)$ with particles as constructed in the proof of Lemma \ref{definition of the map}.
    \end{remark}

    \subsection{The injectivity of the map $\phi_{K}$.} We prove this property by applying the Maximum Principle outside the singular locus and the specialized analysis near cone singularities.

    \begin{proposition}\label{principal curvature at cone singularity}
         Let $(N,g)\in\mathcal{GH}'_{\Sigma,\theta}$ be a convex GHM AdS manifold with particles. Assume that $S$ is a future-convex, spacelike, constant curvature $K$ surface which is orthogonal to the singular lines. Then for each intersection point $p_{i}$ of the surface $S$ with the singular line $l_{i}$ in $N$, both principal curvatures on $S$ tend to $k=\sqrt{-1-K}$ at $p_{i}$ for $i=1,...,n_0$.
    \end{proposition}

    \begin{proof}
        Let $I$ and $B$ be the induced metric and the shape operator of $S$ in $(N,g)$, respectively. Then we have
        \begin{equation*}
            I=\frac{1}{|K|}h,\qquad III=\frac{1}{|K^{*}|}h',
        \end{equation*}
        where $h,h'\in\mathfrak{M}^{\theta}_{-1}$ and $K^{*}=-K/(1+K)$.

        We claim that $id:(S,h)\rightarrow (S,h')$ is minimal Lagrangian. Indeed, set
        \begin{equation*}
            b=\frac{1}{\sqrt{-1-K}}B.
        \end{equation*}

        Note that $S$ is future-convex, then $B$ is positive definite. One can easily check that
        \begin{itemize}
            \item $b$ is self-adjoint for $h$ with positive eigenvalues.
            \item $\det(b)=1$.
            \item $d^{\nabla}b=0$, where $\nabla$ is the Levi-Civita connection of $h$.
            \item $h'=h(b\bullet,b\bullet)$.
       \end{itemize}

       Moreover, by Proposition \ref{Toulisse2}, both eigenvalues of $b$ tend to $1$ at cone singularities. Hence, both eigenvalues of $B=\sqrt{-1-K}b$ tend to $k=\sqrt{-1-K}$ at the intersections of $S$ with the singular lines. This implies the conclusion.
\end{proof}

    Let $(N,g)$ be a convex GHM AdS manifold with particles. Recall that a \emph{convex subset} $\Omega$ of $N$ is a subset of $N$ such that any geodesic segment in $N$ with endpoints in $\Omega$ is contained in $N$. It is proved in \cite[Lemma 4.5, Lemma 4.9]{BS} that the following properties still hold for the case of convex GHM AdS manifolds with particles.

    \begin{lemma}
        \label{the convex core for the case with particles}
            Each convex GHM AdS manifold $(N,g)$ with particles contains a convex core $C(N)$, which is the non-empty convex subset of $N$ such that any non-empty convex subset $\Omega$ in $N$ contains $C(N)$. Moreover, for any point $x\in N\setminus C(N)$, the maximal timelike geodesic segment connecting $x$ to $C(N)$ has length less than $\pi/2$.
    \end{lemma}

    \begin{remark}
        \label{the property of convex core}
The boundary of $C(N)$ is the union of two (possibly identified) surfaces, 
called the future boundary $\partial_{+}C(N)$ and the past boundary $\partial_{-}C(N)$.
            In the Fuchsian case (i.e. the two metrics of the Mess parametrization are equal), $C(N)=\partial_{+}(N)=\partial_{-}(N)$ is a totally geodesic spacelike surface orthogonal to the singular lines. In the non-Fuchsian case,  each boundary component of $C(N)$ is a spacelike surface orthogonal to the singular lines and is ``pleated'' along a measured geodesic lamination. In both cases, the induced metric on each boundary component of $C(N)$ is hyperbolic, with each cone singularity of angle equal to that of corresponding particle, as in \cite[Lemma 1.5]{BS}. Moreover, the maximal geodesic segment starting from $x\in\partial_{+}C(N)$ (resp. $\partial_{-}C(N)$) in the direction of a past-oriented (resp. future-directed) normal vector at $x$ has length $\pi/2$, see \cite[Lemma 1.6]{BS}.
    \end{remark}

    The following corollary is an immediate consequence of Lemma \ref{the convex core for the case with particles}.

    \begin{corollary}
        \label{corollary of the lemma for convex core}
            Let $(N,g)$ be a convex GHM AdS manifold with particles.
            \begin{enumerate}[(1)]
                \item  If $S$ is a strictly future-convex spacelike surface orthogonal to the singular lines in $N$, then $S$ is in the past of the convex core $C(N)$ and stays at distance less than $\pi/2$ from $\partial_{+}C(N)$.
                \item  If $S$ be a strictly past-convex spacelike surface orthogonal to the singular lines in $N$, then $S$ is in the future of the convex core $C(N)$ and stays at distance less than $\pi/2$ from $\partial_{-}C(N)$.
            \end{enumerate}
    \end{corollary}

    The following theorem is an alternative version of the Maximum Principle Theorem (see \cite[Lemma 2.3]{BBZ1}, \cite[Proposition 4.6]{BBZ2}) for the case of convex GHM AdS manifolds with particles.

    \begin{theorem}(Maximum Principle)
        \label{Maximum principle}
            Let $(N,g)$ be a convex GHM AdS manifold with particles. Let $S$ and $S'$ be two future-convex spacelike surfaces in $N$ which are orthogonal to the singular lines. Assume that $S$ and $S'$ intersect at a point $p$ which is not a singularity, and assume that $S'$ is contained in the future of $S$. Then the principal curvatures of $S'$ at $p$ are larger than or equal to those of $S$.
    \end{theorem}

    \begin{lemma}
        \label{the property of the pushing surface}
            Let $(N,g)$ be a convex GHM AdS manifold with particles. Let $S$ be a future-convex, spacelike surface in $N$ orthogonal to the singular lines and let $\psi^{t}:S\rightarrow N$ be a map defined by $\psi^{t}(x)=\exp_{x}(t\cdot n(x))$, where $n(x)$ is the future-directed unit normal vector at $x$ of $S$ in $N$. Then for each $x\in S$ which is a regular point, we have

            \begin{enumerate}[(1)]
                \item $\psi^{t}$ is an embedding in a neighbourhood of $x$ if $t$ satisfies that $\lambda(x)\tan(t)\not=-1$ and $\mu(x)\tan(t)\not=-1$, where $\lambda(x)$ and $\mu(x)$ are the principal curvatures of $S$ at $x$.
                \item  The principal curvatures of $\psi^{t}(S)$ at the point $\psi^{t}(x)$ are given by
                    \begin{equation*}
                        \lambda^{t}(\psi^{t}(x))
                        =\frac{\lambda(x)-\tan (t)}{1+\lambda(x)\tan (t)},\qquad
                        \mu^{t}(\psi^{t}(x))=\frac{\mu(x)-\tan (t)}{1+\mu(x)\tan (t)}.
                    \end{equation*}
                \item Fix $x\in S$, $\lambda^{t}(\psi^{t}(x))$ and $\mu^{t}(\psi^{t}(x))$ are both strictly decreasing in $t\in(t_{0}(x)-\pi/2,\pi/2)$, where $t_{0}(x)=\min\{\arctan\lambda(x),\arctan\mu(x)\}$.
            \end{enumerate}
    \end{lemma}

    \begin{lemma}
         \label{comparison between principal curvatures at singularities}
            Let $(N,g)$ be a convex GHM AdS manifold with particles. Let $S_{1},S_{2}$ be two spacelike surfaces in $N$ which are orthogonal to the singular lines. Assume that $S$ and $S'$ intersect at a singular point $p$ such that the limits of both principal curvatures of $S$ at $p$ are equal to $k>0$, and the limits of both principal curvatures of $S'$ at $p$ are equal to $k'>0$. If there exists a neighbourhood $U$ of $p$ in $S$ and a neighbourhood $U'$ of $p$ in $S'$ such that $U'$ is in the future of $U$, then $k'\geq k$.
    \end{lemma}

    \begin{lemma}\label{curvatures and positions of surfaces}
        Let $S_{i}$ be a future-convex spacelike surface of constant curvature $K_{i}$ which is orthogonal to the singular lines in a convex GHM AdS manifold with particles for $i=1,2$. Then we have the following statements:
        \begin{enumerate}[(1)]
            \item $K_{1}<K_{2}$ if and only if $S_{1}$ is strictly in the past of $S_{2}$.
            \item $K_{1}>K_{2}$ if and only if $S_{1}$ is strictly in the future of $S_{2}$.
            \item $K_{1}=K_{2}$ if and only if $S_{1}$ coincides $S_{2}$.
        \end{enumerate}
    \end{lemma}

    \begin{proof}
        By the symmetry between Statement (1) and Statement (2), it suffices to prove Statement (1).

        First we prove that $K_{1}<K_{2}$ implies that $S_{1}$ is strictly in the past of $S_{2}$. We argue by contradiction. Assume that $S_{1}$ is not strictly in the past of $S_{2}$. Set $t_{0}=\sup\{d(x, S_{1}):x\in S_{2}$ is in the past of $S_{1}\}$, where $d(x,S_{1})$ is the maximum of the Lorentzian lengths of causal segments connecting $x$ to $S_{1}$. It is clear that $t_{0}>0$.

        Note that $d(x,S_{1})$ is continuous (see Lemma 4.3 in \cite{BS}) and $S_{1}$ is compact, thus $t_{0}$ is attained at some point $x_{0}\in S_{2}$. In particular, if $x_{0}$ is a regular point, the distance $t_{0}$ is realized by a geodesic segment with the endpoints orthogonal to $S_{1}$ and $S_{2}$ which avoids the singularities. If $x_{0}$ is a singular point, the distance $t_{0}$ is realized by the segment contained in the singular line through $x_{0}$ which connects $x_{0}$ to $S_{1}$.

        Denote $S_{2}^{t}=\psi^{t}(S_{2})$, where $\psi^{t}$ is the map defined in Lemma \ref{the property of the pushing surface}. Consider $S_{2}^{t_{0}}$, it intersects $S_{1}$ the point $y_{0}=\psi^{t_0}(x_{0})$ and it is in the future of $S_{1}$. We discuss it in the following two cases.

         \vspace{1mm}
         \textbf{Case 1:} $x_{0}$ is a regular point. By Corollary
         \ref{corollary of the lemma for convex core}, $t_{0}\in (0,\pi/2)$. By Statement (3) of
         Lemma \ref{the property of the pushing surface}, $\lambda_{2}^{t}(\psi^{t}(x_{0}))$ and $\mu_{2}^{t}(\psi^{t}(x_{0}))$ are both strictly decreasing in $t\in(0,\pi/2)$. Then we have
         \begin{equation}
            \label{eq1:principal curvatures}
                \lambda_{2}^{t_{0}}(y_0)\mu_{2}^{t_{0}}(y_0)
                <\lambda_{2}(x_{0})
                \mu_{2}(x_{0})=-1-K_{2}.
         \end{equation}

         On the other hand, Theorem \ref{Maximum principle} implies that
         \begin{equation}
            \label{eq2:principal curvatures}
                \lambda_{2}^{t_{0}}(y_0)\mu_{2}^{t_{0}}(y_0)
                \geq \lambda_{1}(y_0)\mu_{1}(y_0)=-1-K_{1},
                \end{equation}
         where $\lambda_{2}^{t_{0}}(y_{0})\geq\lambda_{1}(y_{0})>0$ and $\mu_{2}^{t_{0}}(y_{0})\geq\mu_{1}(y_{0})>0$. Combining \eqref{eq1:principal curvatures} and \eqref{eq2:principal curvatures}, we get $K_{1}> K_{2}$, which contradicts the assumption.

         \vspace{1mm}
         \textbf{Case 2:} $x_{0}$ is a singularity. Note that $S_{2}^{t}$ is obtained by pushing $S_{2}$ along the orthogonal geodesics in the future direction.  In particular, the singularities on $S_{2}^{t}$ are the image of the singularities on $S_{2}$ by pushing along the singular lines.
         By Proposition \ref{principal curvature at cone singularity} and
         Lemma \ref{the property of the pushing surface}, the limits of both principal curvatures of $S_{2}$ at the singularity $x_0$ are equal to
         \begin{equation*}
            \lambda_{2}(x_{0})
            =\mu_{2}(x_{0}):=k_{2},
         \end{equation*}
         and the limits of both principal curvatures of $S_{2}^{t_{0}}$ at $y_0=\psi^{t_0}(x_0)$ are equal to
         \begin{equation}
            \label{ineq1:principal curvature at singularity}
                \lambda_{2}^{t_{0}}(y_{0})
                =\mu_{2}^{t_{0}}(y_{0}):=k_{2}^{t_{0}}< k_{2}
         \end{equation}
         where $k_{2}=\sqrt{-1-K_{2}}$ and $k_{2}^{t_{0}}=(k_{2}-\tan(t_{0}))/(1+k_{2}\tan(t_{0}))$.
         Moreover, $S_{1}$ and $S_{2}^{t_{0}}$ intersects at a singularity $y_{0}$ and $S_{2}^{t_{0}}$ is in the future of $S_{1}$. By Lemma
         \ref{comparison between principal curvatures at singularities}, we have
         \begin{equation}
            \label{ineq2:principal curvature at singularity}
                k_{2}^{t_{0}}\geq k_{1}=\sqrt{-1-K_{1}}.
         \end{equation}
         Combining \eqref{ineq1:principal curvature at singularity} and \eqref{ineq2:principal curvature at singularity}, we have $K_{1}>K_{2}$, which leads to a contradiction.

         Now we prove the sufficiency, that is, if $S_{1}$ is strictly in the past of $S_{2}$, then $K_{1}<K_{2}$. Denote $S_{1}^{t}=\psi^{t}(S_{1})$. Set
         $\delta_{0}=\sup\{d(z,S_{2}): z\in S_{1}\}$. Obviously, $\delta_{0}>0$. Assume $\delta_0$ is attained at a point $z_0\in S_1$ and denote $w_0=\psi^{\delta_0}(z_0)\in S_2\cap S_1^{\delta_0}$. Using the similar argument again, we have
         \begin{equation*}
            \begin{split}
                &\lambda_{1}^{\delta_{0}}(w_{0})\mu_{1}^{\delta_{0}}(w_{0})
                <\lambda_{1}(z_{0})\mu_{1}(z_{0})
                =-1-K_{1},\\
               &\lambda_{1}^{\delta_{0}}(w_{0})\mu_{1}^{\delta_{0}}(w_{0})\geq \lambda_{2}(w_{0})\mu_{2}(w_{0})
               =-1-K_{2}.
            \end{split}
         \end{equation*}
         Thus $K_{1}<K_{2}$.

        Now we prove Statement (3). The sufficiency is obvious. Now we show the necessity.

        By assumption, $K_{1}=K_{2}$. Set $d_{1}=\sup\{d(x,S_{1}):x\in S_{2}$ is in the past of $S_{1}\}$ and $d_{2}=\sup\{d(x,S_{1}):x\in S_{2}$ is in the future of $S_{1}\}$. Note that $S_{1}=S_{2}$ if and only if $d_{1}=d_{2}=0$.

        If $d_{1}>0$, consider the surface $S_{2}^{d_{1}}$ obtained by pushing $S_{2}$ along orthogonal geodesics in a distance $d_{1}$ in the future direction. Using the argument as above, we obtain the contradiction that $K_{2}<K_{1}$. This implies that $d_{1}=0$.

        If $d_{2}>0$,  we consider the surface $S_{1}^{d_{2}}$ obtained by pushing $S_{1}$ along orthogonal geodesics in a distance $d_{2}$ in the future direction. Using the same argument as above, we obtain the contradiction that $K_{1}<K_{2}$. This implies that $d_{2}=0$. Therefore, $S_{1}=S_{2}$.
    \end{proof}

    \begin{lemma}\label{injectivity of the map}
        For any $K\in(-\infty, -1)$, the map $\phi_{K}:\mathcal{T}_{\Sigma,\theta}\times \mathcal{T}_{\Sigma,\theta}\rightarrow \mathcal{GH}_{\Sigma,\theta}$ is injective.
    \end{lemma}

    \begin{proof}
         Assume that $(h,h')$, $(h_{1},h_{1}')\in\mathcal{T}_{\Sigma,\theta}\times \mathcal{T}_{\Sigma,\theta}$ satisfy that $\phi_{K}(h,h')=\phi_{K}(h_{1},h_{1}'):=(N,g)$. Then $(N,g)$ contains a future-convex, spacelike surface $S$ of constant curvature $K$ orthogonal to the singular lines, with the induced metric $I=(1/|K|)h$ and the third fundamental form $III=(1/|K^{*}|)h'$ and a future-convex, spacelike surface $S_{1}$ of constant curvature $K$ orthogonal to the singular lines, with the induced metric $I_{1}=(1/|K|)h_{1}$ and the third fundamental form $III=(1/|K^{*}|)h_{1}'$.
         By Lemma \ref{curvatures and positions of surfaces}, we have $S=S_{1}$. Then $h=h_{1}$ and $h'=h_{1}'$, which implies that $(h,h')=(h_{1},h_{1}')$.
    \end{proof}

    \subsection{The continuity of the map $\phi_{K}$.}

    To see this, we relate minimal Lagrangian maps to harmonic maps and use some basic facts on the properties of harmonic maps and energy.

    Let $f:(M,g)\rightarrow(N,h)$ be a $C^1$ map between two closed Riemannian surfaces (possibly with punctures). The differential $df$ of $f$ is a section of $T^{*}M\otimes f^{*}(TN)$ with the metric $g^{*}\otimes f^{*}h$, where $g^{*}$ is the metric on $T^{*}M$ dual to $g$. The \emph{energy} of $f$ is defined as
    \begin{equation*}
        E(f,g,h)=\int_{M} e(f)\,d\sigma_g,
    \end{equation*}
    where $d\sigma_g$ is the area element of $(M,g)$, and $e(f)=\frac{1}{2}||df||^2_{g^{*}\otimes f^{*}h}$ is called the \emph{energy density} of $f$. We call $f$ a \emph{harmonic} map if it is a critical point of the energy $E$.

    It is known that the value of the energy functional $E$ at such a triple $(f,g,h)$ depends only on the conformal class of g. In particular, set $M=N=\Sigma$ and $g,h\in\mathfrak{M}^{\theta}_{-1}$, the energy functional $E$ depends only on the conformal class $c$ of $g$ (see \cite[equality (3.4)]{GR}). This implies that the harmonicity is conformally invariant on the domain.

    The \emph{Hopf differential} of $f$ is defined as the (2,0) part of the pull-back by $f$ of $h$ in the conformal coordinate of $c$, which is denoted by $\Phi(f)$. It measures the difference between the conformal class of $f^{*}(h)$ and $c$. It is shown (cf. \cite[Lemma 5.1]{GR}) that for $f$ harmonic, $\Phi(f)$ is a meromorphic quadratic differential on $(\Sigma,c)$ with at most simple poles at cone singularities.

    \begin{theorem}\emph{(J. Gell-Redman \cite[Theorem 2]{GR})}
        \label{existence and uniqueness of harmonic maps}
            Given $g\in\mathfrak{M}^{\theta}_{-1}$ and $c\in\mathcal{T}_{\Sigma,\theta}$, there exists a unique harmonic map $u_{c,g}:(\Sigma,c)\rightarrow (\Sigma,g)$ isotopic to identity fixing each marked point, and $u_{c,g}$ is a diffeomorphism on $\Sigma_{\mathfrak{p}}$. Moreover, the harmonic maps $u_{c,g}$ vary smoothly with respect to the target metric $g$.
    \end{theorem}

    Minimal Lagrangian maps between hyperbolic surfaces (with cone singularities of angles less than $\pi$) are related to harmonic maps (see e.g. \cite{Toulisse1,BMS1,Schoen}).

    \begin{theorem}\emph{(Toulisse \cite[Theorem 6.4]{Toulisse1})}
         \label{minnimal Lagrangian and harmonic maps}
        Let $h_1, h_2\in\mathfrak{M}^{\theta}_{-1}$. Then there exists a unique conformal structure $c$ on $\Sigma$ such that
            \begin{equation*}
                \Phi(u_{1})+\Phi(u_{2})=0,
            \end{equation*}
            where $\Phi(u_{i})$ is the Hopf differential of the unique harmonic map $u_i:(\Sigma,c)\rightarrow (\Sigma,h_{i})$ isotopic to the identity for $i=1,2$. Moreover, the map $u_2\circ u_1^{-1}:(\Sigma,h_1)\rightarrow (\Sigma,h_2)$ is minimal Lagrangian and isotopic to the identity.
    \end{theorem}

    Fix $g_0\in\mathfrak{M}^{\theta}_{-1}$. Let $E(\bullet,g_0):\mathcal{T}_{\Sigma,\theta}\rightarrow \mathbb{R}$ be a map which assigns to $c\in\mathcal{T}_{\Sigma,\theta}$ the energy of the (unique) harmonic map $u_{c,g_0}$ as indicated in \ref{existence and uniqueness of harmonic maps}.

    It is known that (see \cite[Proposition 2.14]{Toulisse2}) for each $c\in\mathcal{T}_{\Sigma,\theta}$, the tangent space $T_{c}\mathcal{T}_{\Sigma,\theta}$ of $\mathcal{T}_{\Sigma,\theta}$ at $c$ consists of those trace free, divergence free symmetric (0,2)-tensors on $\Sigma_{\mathfrak{p}}$ of class $\mathcal{H}^2$ and $\mathcal{C}^{2}$. It is identified with the space $\mathcal{QD}_{c}(\Sigma)$ of meromorphic quadratic differentials (with respect to the complex structure $c$) on $\Sigma$ with at most simple poles at singularities, by assigning $q\in\mathcal{QD}_{c}(\Sigma)$ to the real part $\Re(q) \in T_{c}\mathcal{T}_{\Sigma,\theta}$.

    Recall that the $L^2$-metric defined on $T_{c}\mathcal{T}_{\Sigma,\theta}$ is given by the inner product:
    \begin{equation*}
        \langle\langle h,k\rangle\rangle_c=\frac{1}{2}\int_M \trace(HK)d\mu_g,
    \end{equation*}
    where $H$, $K$ are the (1,1)-tensors on $\Sigma_{\mathfrak{p}}$ obtained from $h$ and $k$ via the representative metric $g$ of $c$ (by raising an index), $\mu_g$ is the volume element induced on $\Sigma_{\mathfrak{p}}$ by $g$.

    Let $\xi dz^2, \eta dz^2\in\mathcal{QD}_{c}(\Sigma)$, where $g=\lambda |dz|^2$ under the conformal coordinate $z=x+iy$. As an analog in Teichm\"uller space of closed surfaces (see Section 2.6 in \cite{Tromba}), the Weil-Petersson metric on $\mathcal{T}_{\Sigma,\theta}$ is defined as
    \begin{equation*}
        \langle\xi,\eta\rangle_{WP}=\Re\int_{\Sigma}\frac{\xi\bar{\eta}}{\lambda}dxdy.
    \end{equation*}
    One can check that the Weil-Petersson metric on $T_{c}\mathcal{T}_{\Sigma,\theta}$ is equal to the $L^2$-metric :
    \begin{equation*}
        \langle\xi,\eta\rangle_{WP}=\langle\langle \Re(\xi),\Re(\eta)\rangle\rangle_c.
    \end{equation*}

    The following lemma provides the properties of $E(\bullet,g_0)$ we need, see \cite[Theorem 3.2]{Toulisse2} and \cite[Theorem 3.1.3]{Tromba}.

    \begin{lemma}\label{properties of energy functional}
        $E(\bullet,g_0)$ has the following properties:
        \begin{enumerate}[(1)]
            \item $E(\bullet,g_0)$ is proper.
            \item The Weil-Petersson gradient $\nabla E(\bullet,g_0)(g)$ of $E(\bullet,g_0)$ at $g\in\mathcal{T}_{\Sigma,\theta}$ is (up to a factor) $\Re(\Phi(u_{g,g_0}))$.
            \item The second derivative of $E(\bullet,g_0)$ at a critical point is (up to a positive factor) Weil-Petersson metric (hence, positive definite).
            \item The isotopy class associated to $g_0$ is the only critical point of $E(\bullet,g_0)$.
        \end{enumerate}
    \end{lemma}

    Let $h_1,h_2\in\mathfrak{M}^{\theta}_{-1}$. Define the functional $E_{h_1,h_2}(\bullet)=E(\bullet,h_1)+E(\bullet,h_2)$ over $\mathcal{T}_{\Sigma,\theta}$. By Lemma \ref{properties of energy functional}, $E_{h_1,h_2}(\bullet)$ is proper and has a unique critical point $c\in\mathcal{T}_{\Sigma,\theta}$ such that $\Phi(u_{c,h_1})+\Phi(u_{c,h_2})=0$. As a consequence, we have the following proposition.

    \begin{proposition}
        \label{minimum point of the energy functional}
            The conformal structure $c$ in Theorem \ref{minnimal Lagrangian and harmonic maps} is the unique critical (minimum) point of the functional $E_{h_1,h_2}(\bullet):\mathcal{T}_{\Sigma,\theta}\rightarrow \mathbb{R}$.
    \end{proposition}

    \begin{lemma}\label{Continuity of the map}
        For any $K\in(-\infty, -1)$, the map $\phi_{K}:\mathcal{T}_{\Sigma,\theta}\times \mathcal{T}_{\Sigma,\theta}\rightarrow \mathcal{GH}_{\Sigma,\theta}$ is continuous.
    \end{lemma}

    \begin{proof}
         It suffices to prove that if the sequence $(h_{k},h_{k}')_{k\in \N}$ converges to $(h,h')\in\mathcal{T}_{\Sigma,\theta}\times \mathcal{T}_{\Sigma,\theta}$, then the sequence $(\phi_{K}(h_{k},h_{k}'))_{k\in \N}$ converges to $\phi_{K}(h,h')\in\mathcal{GH}_{\Sigma,\theta}$. Denote by $m_{k}$ the unique minimal Lagrangian map between $(\Sigma, h_{k})$ and $(\Sigma, h_{k}')$ isotopic to the identity and by $m$ the unique minimal Lagrangian map between $(\Sigma, h)$ and $(\Sigma, h')$ isotopic to the identity.

        We claim that the sequence $(m_{k})_{k\in \N}$ converges $m$. Indeed, by Proposition \ref{minimum point of the energy functional}, denote by $c_k$ the unique critical point of $E_{h_k,h'_k}(\bullet)$ and by $c$ the unique critical point of $E_{h,h'}(\bullet)$.

        Now we prove that $c_{k}$ converges to $c$. By Theorem \ref{existence and uniqueness of harmonic maps}, $u_{c,h_k}$ (resp. $u_{c,h'_k}$) vary smoothly with respect to the target metrics $h_k$ (resp. $h'_k$). Combined with $(h_{k},h_{k}')_{k\in \N}\to(h,h')$, it follows that $(E_{h_{k},h_{k}'}(\bullet))_{k\in \N}$ converges to $E_{h,h'}(\bullet)$ in the $C^{1}$ sense, that is,
        \begin{equation*}
            E_{h_{k},h_{k}'}(\bullet)\to E_{h,h'}(\bullet)
        \end{equation*}
         pointwise as $k\to \infty$, and
         \begin{equation*}
            \nabla E_{h_k,h'_k}(\bullet)=C\Re(\Phi(u_{\bullet,h_k})+\Phi(u_{\bullet,h'_k}))\to C\Re(\Phi(u_{\bullet,h})+\Phi(u_{\bullet,h'}))=\nabla E_{h,h'}(\bullet)
         \end{equation*}
         pointwise as $k\to \infty$, where $C$ is a non-zero constant. Note that $E_{h,h'}(\bullet)$ has non-degenerate second derivative at $c$. By the Implicit function theorem on Banach spaces (see \cite[Theorem 2.5.7]{AMT}), $c$ is the limit of the critical points $c_{k}$ of $E_{h_{k},h_{k}'}(\bullet)$. By a closeness result for harmonic maps (see Theorem 7.1 in \cite{GR}), $u_{c_k,h_k}$ (resp. $u_{c_k,h'_k}$) converges to the harmonic map $u_{c,h}$ (resp. $u_{c,h'}$). Combined with Theorem \ref{minnimal Lagrangian and harmonic maps}, $m_k=u_{c_k,h_k}\circ (u_{c_k,h'_k})^{-1}$ converges to $u_{c,h}\circ (u_{c,h'})^{-1}=m$.

        Let $b_{k}:T\Sigma\rightarrow T\Sigma$ be the bundle morphism defined outside the singular locus which is described in Proposition \ref{Toulisse2} with the property that $m_{k}^{*}(h_{k}')=h_{k}(b_{k}\bullet,b_{k}\bullet)$. Then $b_{k}$ converges to a bundle morphism from $T\Sigma$ to $T\Sigma$, which is denoted by $b$.

        Let $I_{k}=(1/|K|)h_k$ and $B_{k}=\sqrt{-1-K}b_k$. Then $(\Sigma, I_{k},B_{k})_{k\in \N}$ converges to $(\Sigma,I,B)$, in the sense that $I_k$, $B_k$ converges to $I=(1/|K|)h$ and $B=\sqrt{-1-K}b$, respectively. This implies that $(\phi_{K}(h_{k},h_{k}'))_{k\in \N}$ converges to $\phi_{K}(h,h')$ in $\mathcal{GH}_{\Sigma,\theta}$. The proof is completed.
    \end{proof}

    \begin{proposition}\label{local homeomorphism of the map}
         For any $K\in(-\infty, -1)$, the map $\phi_{K}:\mathcal{T}_{\Sigma,\theta}\times \mathcal{T}_{\Sigma,\theta}\rightarrow \mathcal{GH}_{\Sigma,\theta}$ is a local homeomorphism.
         \end{proposition}

    \begin{proof}
         By the extension of Mess parametrization (see \cite[Theorem 1.4]{BS}),
         $ \mathcal{GH}_{\Sigma,\theta}$ is homeomorphic to $\mathcal{T}_{\Sigma,\theta}\times \mathcal{T}_{\Sigma,\theta}$.
          Thus, $\mathcal{T}_{\Sigma,\theta}\times \mathcal{T}_{\Sigma,\theta}$ and $ \mathcal{GH}_{\Sigma,\theta}$ have the same dimension and have no boundary. Moreover, it follows from
         Lemma \ref{injectivity of the map} and
         Lemma \ref{Continuity of the map} that $\phi_{K}$ is injective and continuous. By the invariance of domain theorem for manifolds, $\phi_{K}$ is a local homeomorphism.
    \end{proof}

    \subsection{The properness of the map $\phi_{K}$.}

    To prove this property of $\phi_K$, we recall some elementary facts about hyperbolic surfaces with cone singularities of angles less than $\pi$.

    First we introduce the following Collar Lemma for hyperbolic cone surfaces (see \cite[Theorem 3]{DP}).

    \begin{lemma}\emph{(Collar lemma)}\label{Collar lemma for hyperbolic cone-surfaces}
        Let $S$ be a hyperbolic cone-surface of genus $g$ with $n_0$ cone points $p_1,..., p_{n_0}$ with cone angles $\theta_1,...,\theta_{n_0}\in(0,\pi)$ and $(g,{n_0})\geq (0,4)$. Let $\alpha$ be the largest cone angle. Let $\{\gamma_1,...,\gamma_{m}\}$ be a maximal collection of mutually disjoint simple closed geodesics on $S$, where $m=3g-3+n_0$. Then the collars
        \begin{equation*}
            C(\gamma_{k})=\{x\in S:d(x,\gamma_{k})\leq w_{k}=\arcsinh\left(\cos\alpha/\sinh{\frac{\ell(\gamma_{k})}{2}}\right)\}
        \end{equation*}
        and
        \begin{equation*}
            C(p_{l})=\{x\in S:d(x,p_{l})\leq v_{l}=\arccosh(1/\sin{\theta_{l}})\}
        \end{equation*}
        are pairwise disjoint for $k=1,...,m$ and $l=1,...,n_0$, where $\ell(\gamma_{k})$ is the length of the geodesic $\gamma_{k}$.
        \end{lemma}

    \begin{lemma}\label{unbounded sequence of points in Teichmuller space}
        Let $(\tau_{i})_{i\in\mathbb{N}}\subset \mathcal{T}_{\Sigma,\theta}$ be a sequence which escape from any compact subset of $\mathcal{T}_{\Sigma,\theta}$. Then there exists a simple closed curve $\gamma$ on $\Sigma$ such that, up to extracting a subsequence, the length of $\gamma$ under $\tau_{i}$ tends to infinity.
    \end{lemma}

    \begin{proof}
         Note that the underlying surface $\Sigma$ we consider satisfies the condition
         \begin{equation*}
            2\pi(2-2g)+\sum\limits_{i=1}\limits^{n}(\theta_{i}-2\pi)<0.
         \end{equation*}
         Then each marked hyperbolic cone-surface in $\mathcal{T}_{\Sigma,\theta}$ admits a pants decomposition $\mathcal{P}=\{C_{i}\}_{i=1}^{3g-3+n_0}$ such that each pair of pants obtained from $\mathcal{P}$ is either a hyperbolic pair of pants with three boundary components or a generalized hyperbolic pair of pants with exactly one or two boundary components degenerating into cone points of the given angles. In the latter case, the pair of pants is uniquely determined by the lengths of the non-degenerated boundary components and the angles of the cone points.

         With the angles of the cone points fixed, $\mathcal{T}_{\Sigma,\theta}$ has the analogous Fenchel-Nielsen coordinate as the usual Teichm\"uller space of hyperbolic surfaces. Moreover, the twist parameter along each pant curve $C_{i}$ is determined by the length of the shortest simple closed geodesic $\alpha_{i}$ which intersects $C_{i}$ and the length of the geodesic $T_{C_{i}}(\alpha_{i})$ obtained by taking a positive Dehn-twist along $C_{i}$ on $\alpha_{i}$. This implies that there exist finitely many simple closed curves on $\Sigma$ whose lengths completely determine an element in $\mathcal{T}_{\Sigma,\theta}$.

         By assumption, $(\tau_{i})_{i\in\mathbb{N}}$ escape from any compact subset of $\mathcal{T}_{\Sigma,\theta}$. Then there must be some simple closed curve $\gamma$ whose length under $\tau_{i}$ tends to either infinite or zero (in the latter case, it follows from Lemma \ref{Collar lemma for hyperbolic cone-surfaces} that any simple closed curve intersecting $\gamma$ is becoming infinitely long). Therefore, there always exists a simple closed curve (still denoted by $\gamma$) on $\Sigma$, such that up to extracting subsequence, the length of $\gamma$ under $\tau_{i}$ tends to infinity. This completes the proof.
    \end{proof}

    The following lemma gives a comparison between the lengths of simple closed geodesics in the same isotopy class on the past boundary $\partial_{-}C(N)$ of the convex core and on a spacelike surface in its past in a convex GHM AdS manifold $(N,g)$.

    \begin{lemma}
        \label{the comaprison of the length on future-convex surfaces}
            Let $(N,g)$ be a convex GHM AdS manifold with particles. Let $S$ be a spacelike surface in the past of $\partial_{-}C(N)$ which is orthogonal to the singular lines in $N$. Then for any closed geodesic $\gamma$ on $\partial_{-}C(N)$, the length of $\gamma$ is larger than the length of any closed minimizing geodesic $\gamma'$ on $S$ homotopic to $\gamma$.
    \end{lemma}

    \begin{proof}

        Let $\lambda_{-}$ be the bending lamination of $\partial_{-}C(N)$. 
The set of isotopy classes of weighted non-trival simple closed curves is dense in the space $\cML_{\Sigma,n_0}$ of measured laminations on $\Sigma_{\mathfrak{p}}$ (see \cite[Proposition 3.1]{BS}). It suffices to consider the case where $\lambda_{-}$ is a disjoint finite union of weighted simple closed geodesics on $\partial_{-}C(N)$. Assume that $\supp\lambda_{-}=\cup_{i=1}^{m}\alpha_{i}$, where $\alpha_i$ is a simple closed geodesic on $\partial_{-}C(N)$ disjoint from $\alpha_{j}$ for $j\not=i$. Then $\partial_{-}C(N)\setminus \supp\lambda_{-}$ is a disjoint finite union of spacelike subsurfaces of $\partial_{-}C(N)$ which are totally geodesic in $N$.

Let $\Sigma_{0}=\partial_{-}C(N)$, and let $h_0$ be the induced metric on $\Sigma_0$.
First we construct a family $(\Sigma_{t})_{t\in [0,\pi/2]}$ of future-convex equidistant surfaces from $\partial_{-}C(N)$ in $I^{-}(\Sigma_{0})$. For each $t\in(0,\pi/2]$, let
$$ \Omega_t=\{x\in I^{-}(\Sigma_{0})~|~ d(x,\Sigma_0)\leq t\}~, $$
and let
$$ \Sigma_t = \partial\Omega_t\cap I^-(\Sigma_0)~. $$
Note that $\Sigma_t$ is a future-convex (non-smooth) spacelike surface orthogonal to the singular lines (see e.g. \cite[Lemma 4.2]{BS}) and $\Sigma_t$ can be disconnected when it is close to the past singularity of $N$ and even empty when $t$ tends to $\pi/2$. It is clear that $\cup_{t\in(0,\pi/2)}\Sigma_t=I^{-}(\Sigma_{0})$.

Let $x\in \Sigma_t$, for some $t\in (0,\pi/2)$, and let $n$ be a unit future-oriented vector orthogonal to a support plane of $\Sigma_t$ at $x$. Let $\gamma_{x,n}$ be the intersection with $I^-(\Sigma_0)\cup\Sigma_0$ of the geodesic starting from $x$ with velocity $n$. Since $\Sigma_t$ is future-convex, the $\gamma_{x,n}$ are disjoint. We define an ``orthonormal projection'' $p_t$ to $\Sigma_t$, sending a point $y\in I^+(\Sigma_t)\cap (I^-(\Sigma_0)\cup\Sigma_0)$ to $x\in \Sigma_t$ if $y\in \gamma_{x,n}$ for a certain time-like unit vector $n$ orthogonal to a support plane of $\Sigma_t$ at $x$. Since $\Sigma_t$ is future-convex, $x$ is then the unique point on $\Sigma_t$ realizing the distance to $x$. Denote by $Dom(p_t)$ the domain of $p_t$, which is a subset of $I^+(\Sigma_t)\cap (I^-(\Sigma_0)\cup\Sigma_0)$.

Let $r,s>0$ and let $y\in Dom(p_{r+s})$. Then $p_{r+s}(y)=p_s(p_r(y))$, because the time-like geodesic segment between $y$ and $p_{r+s}(y)$ must intersect $\Sigma_r$ at a point which maximizes both the distance between $y$ and $\Sigma_r$ and the distance between $p_r(y)$ and $p_{r+s}(y)$.

It follows that there exists a flow $(\phi_t)_{t\in [0,\pi/2]}$, defined for each $t$ on a subset of $I^-(\Sigma_0)$, such that if $y\in \Sigma_r\cap Dom(p_{r+s})$, then $p_{r+s}(y)=\phi_s(y)$. By definition, $(\phi_t)_{t\in [0,\pi/2]}$ is the flow of a past-oriented unit time-like vector field $X$, which is however not continuous. At each point $x\in \Sigma_r$, for $r\in [0,\pi/2)$, $X$ is normal to a support plane of $\Sigma_r$. Although $X$ is discontinuous, it follows from its definition that the flow of $X$ exists (but the flow of $-X$ is not well-defined).

A direct examination shows that the restriction of $p_r$ to $\Sigma_0$ is distance-decreasing. In fact, regions near a pleating line of $\Sigma_0$ are typically sent to a line, and the length along pleating geodesics is contracted by a factor $\cos(r)$. Similarly, on flat regions of $\Sigma_0$ which are sent to smooth regions on $\Sigma_r$, lengths are contracted by a factor $\cos(r)$. So, if we denote by $h_r $ the pull-back on $\Sigma_0$ by $p_r$ of the induced metric on $\Sigma_r$, then $(h_r)_{r\in (0,\pi/2)}$ is a decreasing family of pseudo-metrics (each defined on a subset of $\Sigma_0$, this subset being 
also decreasing with $r$).

We now consider the map $\phi:\Sigma_0\to S$, with $\phi(x)$ defined by following the flow of $X$ from $x$ to the first intersection point with $S$. For all $x\in \Sigma_0$, we also denote by $t(x)$ the time needed to reach $\phi(x)$, so that $\phi(x)\in \Sigma_{t(x)}$. Finally we denote by $h$ the pseudo-metric obtained on $\Sigma_0$ as the pull-back by $\phi$ of the induced metric on $S$. (Note that $h$ is defined on the whole of $\Sigma_0$
because $\phi$ is defined on the whole of $\Sigma_0$ since any integral curve of $X$ starting from $\Sigma_0$ must intersect $S$. For the same reason, $h_{t(x)}$ is well-defined at $x$.)

Let $x\in \Sigma_0$. At $\phi(x)$, the tangent plane $T_{\phi(x)}S$ can be identified to the tangent $P$ to any support plane of $\Sigma_{t(x)}$ by projection along the normal to $P$. Under this identification, the induced metric on $T_{\phi(x)}S$ is smaller than the induced metric on $P$ (the difference being $dt^2$, where $t$ denotes now the distance to $\Sigma_0$).

It follows that, at all $x\in \Sigma_0$, $h\leq h_{t(x)}$, and therefore $h\leq h_0$.

Let $\gamma$ be a closed geodesic on $\Sigma_0$, and let $\gamma'=\phi(\gamma)\subset S$. The length of $\gamma$ for $h$ is smaller than the length of $\gamma$ for $h_0$, so that the length of $\gamma'$ for the induced metric on $S$ is less than the length of $\gamma$ for the induced metric on $\Sigma_0$. It follows that the length on $S$ of any minimizing geodesic homotopic to $\gamma'$ (and therefore to $\gamma$) is smaller than the length of $\gamma$.
\end{proof}

    Note that a much simpler proof of the 
    Lemma \ref{the comaprison of the length on future-convex surfaces} can be given if $S$ is a future-convex spacelike surface orthogonal to the singular lines and if there is a foliation of the region between $\partial_-C(N)$ and 
    $S$ by smooth (outside the singular locus) future-convex surfaces orthogonal to the singular lines. The existence of such a foliation clearly follows from Theorem \ref{foliation of the GHMC AdS manifold}. However, 
    at this point of the proof, we couldn't find a simple way to prove the existence of such a foliation by smooth future-convex surfaces. Therefore, we give an alternative method, which also generalizes the case of a future-convex surface $S$ (orthogonal to the singular lines) to the case of a spacelike surface (orthogonal to the singular lines) in the past of $\partial_{-} C(N)$.

    The following corollary is an analogue of Lemma \ref{the comaprison of the length on future-convex surfaces}.

    \begin{corollary}
         \label{the comaprison of the length on past-convex surfaces}
            Let $(N,g)$ be a convex GHM AdS manifold with particles. Let $S$ be a spacelike Cauchy surface in the future of $\partial_{+}C(N)$ which is orthogonal to the singular lines in $N$. Then for any closed geodesic $\gamma$ on $\partial_{+}C(N)$, the length of $\gamma$ is larger than the length of the closed geodesic $\gamma'$ on $S$ homotopic to $\gamma$.
    \end{corollary}

    \begin{proposition}\label{properness of the map}
         For any $K\in(-\infty, -1)$, the map $\phi_{K}:\mathcal{T}_{\Sigma,\theta}\times \mathcal{T}_{\Sigma,\theta}\rightarrow  \mathcal{GH}_{\Sigma,\theta}$ is proper.
    \end{proposition}

    \begin{proof}
        Let $(N_{k},g_{k}):=\phi_{K}(h_{k},h_{k}')$. It suffices to verify that if a sequence $(h_{k},h_{k}')_{k\in \N}$ escape from any compact subset of $\mathcal{T}_{\Sigma,\theta}\times \mathcal{T}_{\Sigma,\theta}$, then $(N_{k},g_{k})_{k\in \N}$ escape from any compact subset of $\mathcal{GH}_{\Sigma,\theta}$. Indeed, if $(h_{k},h_{k}')_{k\in \N}$ escape from any compact subset of $\mathcal{T}_{\Sigma,\theta}\times \mathcal{T}_{\Sigma,\theta}$, then $(h_{k})_{k\in \N}$ or $(h_{k}')_{k\in \N}$ escape from any compact subset of $\mathcal{T}_{\Sigma,\theta}$. We discuss in the following two cases.

       \vspace{1mm}
       \textbf{Case 1:}
       If $(h_{k})_{k\in \N}$ escape from any compact subset of $\mathcal{T}_{\Sigma,\theta}$.  By Lemma \ref{unbounded sequence of points in Teichmuller space}, there is a simple closed curve $\gamma$ on $\Sigma$, such that up to extracting a subsequence, $\ell_{h_{k}}(\gamma)\rightarrow \infty$.

       Denote by $S_{k}$ the future-convex, constant curvature $K$ surface which is orthogonal to the singular lines with the induced metric $I_{k}=(1/|K|)h_{k}$ in $(N_{k},g_{k})$. Denote the induced metric on $\partial_{-}C(N_{k})$ by $I^{-}_{k}$. It is shown in \cite[Lemma 5.4]{BS} that $I^{-}_{k}$ is a hyperbolic metric with cone singularities of angles equal to the given angles at the intersections with the corresponding singular lines. By Lemma \ref{the comaprison of the length on future-convex surfaces}, $\ell_{I^{-}_{k}}(\gamma)\geq\ell_{I_{k}}(\gamma)\rightarrow \infty$. Note that $\mathcal{GH}_{\Sigma,\theta}$ can be parameterized by the embedding data (including the induced metric and the bending lamination) of the past (or future) boundary of the convex core (see e.g. \cite{BS}). This implies that $(N_{k},g_{k})_{k\in \N}$ are not contained in any compact subset of $ \mathcal{GH}_{\Sigma,\theta}$.

       \vspace{1mm}
       \textbf{Case 2:}
       If $(h'_{k})_{k\in \N}$ escape from any compact subset of $\mathcal{T}_{\Sigma,\theta}$. By Lemma \ref{definition of the map}, the future-convex constant curvature $K$ surface $S_{k}$ in $(N_{k},g_{k})$ has third fundamental form $III_{k}=(1/|K^{*}|)h'_{k}$, where $K^{*}=-K/(1+K)$. By Proposition \ref{the duality between convex surfaces}, the dual surface $S^{*}_{k}$ of $S_{k}$ is a past-convex constant curvature $K^{*}$ surface which is orthogonal to the singular lines with the induced metric $I^{*}_{k}$ such that the pull back of $I^{*}_{k}$ on $S^{*}$ through the duality map is $III_{k}$.

       Using a similar argument as in the first case and applying Corollary \ref{the comaprison of the length on past-convex surfaces}, there exists a simple closed curve $\gamma'$ on $\Sigma$, such that up to extracting a subsequence, $\ell_{I^{+}_{k}}(\gamma')\geq\ell_{I^{*}_{k}}(\gamma')\rightarrow \infty$, where $I^{+}_{k}$ denotes the induced metric on $\partial_{+}C(N_{k})$. This implies that $(N_{k},g_{k})_{k\in \N}$ are not contained in any compact subset of $ \mathcal{GH}_{\Sigma,\theta}$.

       Combining these two cases, the proof is complete.
    \end{proof}

    \emph{Proof of Theorem \ref{parametrization map}.}
        Note that $\mathcal{T}_{\Sigma,\theta}\times \mathcal{T}_{\Sigma,\theta}$ and $ \mathcal{GH}_{\Sigma,\theta}$ are simply connected. By Proposition \ref{local homeomorphism of the map} and Proposition \ref{properness of the map}, for each $K<-1$, $\phi_{K}$ is both a local homeomorphism and a proper map, which implies that $\phi_{K}$ is a homeomorphism.

    \section{The existence and uniqueness of foliations.}

    In this section, we prove Theorem \ref{foliation of the GHMC AdS manifold}, as an application of Theorem \ref{parametrization map}. Let $(N,g)$ be a convex GHM AdS manifold with particles. Denote by $B^{+}$ and $B^{-}$ the future and the past component of $N\setminus C(N)$.

    To prove Theorem \ref{foliation of the GHMC AdS manifold}, we first show that $B^{-}$ admits a unique foliation by future-convex constant curvature surfaces orthogonal to the singular lines. Note that there is a duality between future-convex and past-convex surfaces orthogonal to the singular lines in GHM AdS manifolds with particles (see Proposition \ref{the duality between convex surfaces}). It is a direct consequence that $B^{+}$ admits a unique foliation by past-convex constant curvature surfaces orthogonal to the singular lines.

    Indeed, Theorem \ref{parametrization map} says that for each $K\in(-\infty,-1)$, the map $\phi_{K}:\mathcal{T}_{\Sigma,\theta}\times \mathcal{T}_{\Sigma,\theta}\rightarrow \mathcal{GH}_{\Sigma,\theta}$ is a homeomorphism. In particular, $\phi_{K}$ is a surjection. This implies that there exists an embedded future-convex spacelike surface $S_{K}$ of constant curvature $K$ which is orthogonal to the singular lines in $N$. Moreover, it follows from the injectivity of $\phi_{K}$ and Corollary \ref{corollary of the lemma for convex core} that this surface $S_{K}$ is unique and contained in $B^{-}$. This implies that the union of $S_{K}$ over all $K\in(-\infty,-1)$ is contained in $B^{-}$. It remains to show that the union of $S_{K}$ over all $K\in(-\infty,-1)$ is exactly $B^{-}$.

    To prove this, we first generalize the notion of the uniformly spacelike (see Definition 3.7 in \cite{BBZ2}) property of a sequence of spacelike surfaces to the case with cone singularities as follows.

    \begin{definition}
        A sequence $(S_{k})_{k\in\mathbb{N}}$ of spacelike surfaces orthogonal to the singular lines in $N$ is said to be \emph{uniformly spacelike}, if for every sequence $(x_{k})_{k\in\mathbb{N}}$ with $x_{k}\in S_{k}$, it falls into exactly one of the following two classes:
        \begin{enumerate}[(1)]
            \item  $x_{k}\in S_{k}$ escapes from any compact subset of $N$.
            \item  Up to extracting a subsequence, the sequence $(x_{k},P_{k})_{k\in\mathbb{N}}$ converges to a limit $(x,P)$, with $x\in N$ and $P$ a totally geodesic spacelike plane through $x$. Here $P_{k}$ is the tangent plane of $S_{k}$ at $x_{k}$ if $x_{k}$ is a regular point, and $P_{k}$ is the totally geodesic plane orthogonal to the singular line through $x_{k}$ if $x_{k}$ is a singular point. For convenience, we call $P_{k}$ the \emph{support plane} of $S_{k}$ at $x_{k}$ whatever $x_{k}$ is regular or not.
        \end{enumerate}
    \end{definition}

    Let $(S_{k})_{k\in\mathbb{N}}$ be a sequence of future-convex spacelike surfaces orthogonal to the singular lines in $N$, such that $S_{k+1}$ is strictly in the past of $S_{k}$ for all $k\in\mathbb{N}$. Denote by $\Omega$ the union of the future $I^{+}(S_{k})$ of $S_{k}$ over all $k\in\mathbb{N}$ and denote by $S_{\infty}=\partial \Omega$ the boundary of $\Omega$.

    Note that after pushing along geodesics orthogonal to a future-convex spacelike surface $S$ (orthogonal to the singular lines) in the future direction for the distance $t\in[0,\pi/2]$, the obtained surface is still orthogonal to the singular lines. In the case $\Omega\not=N$, the property of $\partial \Omega$ and the uniformly spacelike property of $(S_{k})_{k\in\mathbb{N}}$  (see Theorem 3.6 and Corollary 3.8 in \cite{BBZ2}) can be directly generalized to the case with cone singularities as follows.

    \begin{lemma}\label{the union is not the whole spacetime}
        Let $\Omega$ and $S_{\infty}$ be the domain and the surface in $N$ as described above. Assume that $\Omega\not=N$, then
          \begin{enumerate}[(1)]
            \item  $S_{\infty}$ is the set of limits in $N$ of sequences $(x_{k})_{k\in\mathbb{N}}$ with $x_{k}\in S_{k}$.
            \item  $S_{\infty}$ is a future-convex spacelike surface which is orthogonal to the singular lines in $N$.
            \item $(S_{k})_{k\in\mathbb{N}}$ is uniformly spacelike.
          \end{enumerate}
    \end{lemma}

    To prove Theorem \ref{foliation of the GHMC AdS manifold}, we need the following compactness result, which is an elementary fact about $\mathcal{T}_{\Sigma,\theta}$. In the case of hyperbolic metrics on closed surfaces, we refer to Lemma 9.4 in \cite{BMS2}.

    \begin{lemma}\label{the compactness of some subset of Teichmuller space}
        Let $C>1$ and $h\in \mathcal{T}_{\Sigma,\theta}$. Let $\mathcal{B}(h)$ be the set consisting of $h'\in\mathcal{T}_{\Sigma,\theta}$ such that for all simple closed curves $\gamma$ on $\Sigma$, $\ell_{\gamma}(h')\leq C\ell_{\gamma}(h)$. Then $\mathcal{B}(h)$ is compact.
    \end{lemma}

    \begin{proof}
        We argue by contradiction. It is clear that $\mathcal{B}(h)$ is closed. Suppose that $\mathcal{B}(h)$ is not compact. Then there exists a sequence $(h_{k}')_{k\in \N}\subset\mathcal{B}(h)$ such that $(h_{k}')_{k\in \N}$ escape from any compact subset of $\mathcal{T}_{\Sigma,\theta}$. By Lemma \ref{unbounded sequence of points in Teichmuller space}, there exists a simple closed curve $\gamma_{0}$ on $\Sigma$, such that up to extracting a subsequence,
        $\ell_{\gamma_{0}}(h_{k}')\rightarrow \infty$,
        as $n\rightarrow \infty$.
        This contradicts the fact that
         \begin{equation*}
            \ell_{\gamma_{0}}(h_{k}')\leq C\ell_{\gamma_{0}}(h)<\infty.
         \end{equation*}
    \end{proof}



    \begin{proposition}
        \label{the convergence of the sequence of convex surfaces}
            Let $(N,g)$ be a convex GHM AdS manifold with particles. Let $(S_{i})_{i\in\mathbb{N}^{+}}$ be a sequence of future-convex spacelike surfaces of constant curvatures $K_{i}$ which are orthogonal to the singular lines in $N$, such that $K_{i+1}<K_{i}$ for all $i\in\mathbb{N}^{+}$. Then the following statements hold.
           \begin{enumerate}[(1)]
           \item  If $K_{i}\rightarrow -\infty$, then the union of $I^{+}(S_{i})$ over $i\in\mathbb{N}^{+}$ is exactly the manifold $N$.
             \item If $K_{i}\rightarrow K$ with $-\infty<K<-1$, then the sequence $(S_{i})_{i\in\mathbb{N}^{+}}$ converges to a future-convex spacelike surface $S_{\infty}$ of constant curvature $K$ (which is orthogonal to the singular lines) in the $\mathcal{C}^2$-topology outside the singular locus.
           \end{enumerate}
    \end{proposition}

    \begin{proof}
        Proof of Statement (1): Assume that the union $\Omega$ of $I^{+}(S_{i})$ over all $i\in\mathbb{N}$ is not $N$.
        By Lemma \ref{the union is not the whole spacetime}, the boundary $S_{\infty}=\partial\Omega$ of $\Omega$ is a future-convex spacelike surface. Moreover, $(S_{i})_{i\in\mathbb{N}^{+}}$ is uniformly spacelike. Therefore, the area of $S_{i}$ does not tend to zero as $i\rightarrow \infty$.
        However, by the Gauss-Bonnet formula for surfaces with cone singularities (see e.g. \cite[Propositon 1]{Troyanov}), the area of $S_{i}$ is equal to
        \begin{equation*}
        \frac{2\pi}{K_{i}}\{\chi(\Sigma)+
        \sum\limits_{i=1}^{n_0}(\frac{\theta_{i}}{2\pi}-1)\},
        \end{equation*}
        where $\Sigma$ is the surface such that $N$ is homeomorphic to $\Sigma\times \mathbb{R}$.
        This implies that the area of $S_{i}$ tends to zero, which produces contradiction.

        \vspace{3mm}
        Proof of Statement (2): Denote by $\Omega$ the union of $I^{+}(S_{i})$ over all $i\in\mathbb{N}^{+}$. First we claim that $\Omega$ is not the whole manifold $N$. To see this, we take a number $K'<K$, it follows from Theorem \ref{parametrization map} and Lemma
        \ref{curvatures and positions of surfaces} that there exists an embedded future-convex spacelike surface $S_{K'}\subset N$ of constant curvature $K'$ (which is orthogonal to the singular lines), such that $S_{K'}$ is strictly in the past of the surfaces $S_{i}$ for all $i\in\mathbb{N}^+$. Hence, the closure of $\Omega$ is contained in the closure of $I^{+}(S_{K'})$. This implies that $\Omega\not=N$.

        Denote by $S_{\infty}=\partial\Omega$ the boundary of $\Omega$.
        By Lemma \ref{the union is not the whole spacetime}, $S_{\infty}$ is a future-convex spacelike surface which is orthogonal to the singular lines in $N$.
        Let $S_{0}=\partial_{-}C(N)$. Then $S_0$ is a spacelike surface of constant curvature $K_0=-1$ which is orthogonal to the singular lines.

        Denote by $g_{i}$, $g_{\infty}$ the metrics induced on $S_{i}$, $S_{\infty}$ by the Lorentzian metric $g$ on $N$ for all $i\in\mathbb{N}$.
        Let $f_{\infty}:\Sigma\rightarrow N$ be an embedding such that $f_{\infty}(\Sigma)=S_{\infty}$, where $\Sigma$ is the surface such that $N$ is homeomorphic to $\Sigma\times \mathbb{R}$. Let $\psi_{i}:S_{\infty}\rightarrow S_i$ be the homeomorphism obtained by the Gauss normal flow. Denote $f_i=\psi_i\circ f_\infty$. Then $f_i:\Sigma\rightarrow N$ is an embedding such that $f_i(\Sigma)=S_{i}$ for all $i\in\mathbb{N}$. It suffices to prove that $f_{i}$ converges to $f_{\infty}$ in the $\mathcal{C}^{2}$-topology outside the singular locus and $S_{\infty}$ has constant curvature $K$.

        Note that all the surfaces $S_{i}$ are orthogonal to the singular lines $l_{k}$ (which are homeomorphic to $\{p_{k}\}\times \mathbb{R}$) and the angle of the singularity on $S_{i}$ at the intersection with $l_{k}$ is $\theta_{k}\in(0,\pi)$ for $k=1,...,n_0$.
        Therefore, the metrics $g_{i}$ can be written as follows:
        \begin{equation*}
        g_{i}=(1/|K_{i}|)\widehat{g_{i}},
        \end{equation*}
        where $\widehat{g_{i}}\in\mathfrak{M}^{\theta}_{-1}$ for all $i\in\mathbb{N}$.

         By Lemma  
         \ref{the comaprison of the length on future-convex surfaces},
         for any simple closed curve $\gamma$ on $\Sigma$,
         we have
         \begin{equation*}
         \ell_{f_{i}(\gamma)}(g_{i})\leq\ell_{f_{0}(\gamma)}(g_{0}),
         \end{equation*}
         for all $i\in\mathbb{N}$.
         Note that $K_{i}$ decreasingly converges to $K\in(-\infty, -1)$. Then
         \begin{equation*}
            \ell_{f_{i}(\gamma)}(\widehat{g_{i}})
         =\ell_{f_{i}(\gamma)}(|K_{i}|g_{i})
         =\sqrt{|K_{i}|}\ell_{f_{i}(\gamma)}(g_{i})
         \leq\sqrt{|K|}\,\ell_{f_{0}(\gamma)}(g_{0})
         =\sqrt{K/K_{0}}\,\ell_{f_{0}(\gamma)}(\widehat{g_{0}}),
         \end{equation*}
         for all $i\in\mathbb{N}$. Here $K/K_{0}>1$.

         Denote by $f_{i}^{*}(\widehat{g_{i}})$ the pull back metric on $\Sigma$
         of $\widehat{g_{i}}$ by $f_{i}$ and still denote by $f_{i}^{*}(\widehat{g_{i}})$ its isotopy class in $\mathcal{T}_{\Sigma,\theta}$ for all $i\in\mathbb{N}$.
         For any simple closed curve $\gamma$ on $\Sigma$, we get
         \begin{equation*}
         \ell_{\gamma}(f_{i}^{*}(\widehat{g_{i}})) \leq\sqrt{K/K_{0}}\,\ell_{\gamma}({f_{0}}^{*}(\widehat{g_{0}})).
         \end{equation*}

         By Lemma \ref{the compactness of some subset of Teichmuller space}, the set $\{f_{i}^{*}(\widehat{g_{i}}):i\in\mathbb{N}\}$ is compact in $\mathcal{T}_{\Sigma,\theta}$. This implies that, after extracting a subsequence, $(f_{i}^{*}(\widehat{g_{i}}))_{i\in \N}$ converges to the metric $\widehat{g}_{\infty}$ in the $\mathcal{C}^{2}$-topology outside the singular locus, where $\widehat{g}_{\infty}\in\mathcal{T}_{\Sigma,\theta}.$

         Moreover, it follows from Lemma
         \ref{the union is not the whole spacetime} that
        for each point $x\in S_{\infty}$, there exists a sequence $(x_{i},P_{i})_{i\in\mathbb{N}}$ (where $x_{i}\in S_{i}$ and $P_{i}$ is the spacelike support plane of $S_{i}$ at $x_{i}$) such that $x$ is the limit of $x_{i}$ and $(P_{i})_{i\in\mathbb{N}}$ converges to a spacelike support plane of $S_{\infty}$ at $x$. Therefore, the sequence of 1-jets $(j^{1}(f_{i}(x))_{i\in \N}$ converges at regular points.
        If $f_{i}$ does not converge to $f_{\infty}$ in the $\mathcal{C}^2$-topology, it follows from the proof of \cite[Theorem 5.6]{Schlenker} that $f_{\infty}(\Sigma)=S_{\infty}$ is a degenerate surface pleated along a geodesic $\gamma$ (which may be a geodesic segment between the singular points) and tangent to a light-like plane somewhere outside $\gamma$. This implies that $S_{\infty}$ admits a light-like support plane somewhere, which contradicts the fact that $S_{\infty}$ admits a spacelike support plane at each point. Moreover, the induced metric on $S_{\infty}$ is
         \begin{equation*}
         g_{\infty}=\frac{1}{|K|}f_{*}(\widehat{g}_{\infty}),
         \end{equation*}
         which has constant curvature $K$, where $f_{*}(\widehat{g}_{\infty})$ is the push-forward by $f$ of $\widehat{g}_{\infty}\in\mathcal{T}_{\Sigma,\theta}$. This completes the proof of Proposition \ref{the convergence of the sequence of convex surfaces}.
    \end{proof}

    Recall that the \emph{cosmological time} of a spacetime $(M,g)$ is the function $\tau:M\rightarrow[0,+\infty]$ associating to $x\in M$ the supremum of the Lorentzian lengths of all past-oriented inextensible causal curves starting from $x$. It is said to be \emph{regular} if $\tau(x)<+\infty$ for all $x\in M$ and for each past-oriented inextensible causal curve $\gamma:[0,+\infty)\to M$, the limit $\tau(\gamma(t))\to 0$ as $t\to +\infty$.

    Replace ``past-oriented" by ``future-oriented" in the definition of the cosmological time $\tau$, we define the \emph{reverse of the cosmological time} $\breve{\tau}$.

    In general, the cosmological time of a spacetime is not regular (e.g. Minkowski space and de Sitter space). In our case, a convex GHM AdS spacetime $(N,g)$ with particles has a regular cosmological time $\tau$. By Remark \ref{the property of convex core}, we have $B^{+}=\{x\in N: \tau(x)>\pi/2\}$ and $B^{-}=\{x\in N: \breve{\tau}(x)>\pi/2\}$.

    \begin{proposition}\label{the filling is full}
        Let $(N,g)$ be a convex GHM AdS manifold with particles. Then $B^{-}$ is exactly the union of the surface $S_{K}$ over $K\in(-\infty, -1)$, where $S_{K}$ is the future-convex spacelike surface of constant curvature $K$ which is orthogonal to the singular lines in $N$.
    \end{proposition}

    \begin{proof}
       Denote by $V$ the union of the surface $S_{K}$ over $K\in(-\infty, -1)$. Moreover, Lemma \ref{curvatures and positions of surfaces} implies that $S_K$ is disjoint from $S_{K'}$ for all $K\not=K'\in(-\infty,-1)$.
       Note that $V$ is contained in $B^{-}$. We only need to prove that $B^{-}$ is contained in $V$.

       Fix a number $K_{1}<-1$. Consider the union $V_{1}$ of the surfaces $S_{K}$ over $K\in(-\infty, K_{1})$. By Proposition
       \ref{the convergence of the sequence of convex surfaces}, we have $V_{1}\cap B^{-}=I^-(S_{K_1})\cap B^{-}$. Let $V_{2}=V\setminus V_1$, that is, the union of the surface $S_{K}$ over $K\in[K_{1},-1)$. It is enough to show that $B^{-}\setminus V_{1}\subset V_{2}$. The argument is similar to that of Claim 11.14 in \cite{BBZ2}. For completeness, we include the proof as follows.

        Denote by $V_{2}^{*}$ the union of the surfaces $S_{K}^{*}$ dual to $S_{K}$ over all $K\in[K_{1},-1)$. By Proposition \ref{the duality between convex surfaces}, the surface $S_{K}^{*}$ is a past-convex spacelike surface in $B^{+}$ of constant curvature $K^{*}=-K/(1+K)$, which is orthogonal to the singular lines in $N$. Observe that $K^{*}\rightarrow \infty$ iff $K\rightarrow -1$.

        Note that Proposition \ref{the convergence of the sequence of convex surfaces} is applicable to the family $\{S_{K^{*}}:=S_{K}^{*}\}_{K^{*}\in (-\infty, K_1^{*}]}\subset B^{+}$ (it follows directly from reversing the time orientation of $N$), where $K_{1}^{*}=-K_{1}/(1+K_1)$. This implies that
        \begin{equation*}
            \lim_{K^{*}\rightarrow-\infty}\sup\limits_{x\in S_{K^{*}}}\breve{\tau}(x)
            =\lim_{K\rightarrow -1}\sup\limits_{x\in S^{*}_{K}}\breve{\tau}(x)
            =0.
        \end{equation*}
         where $\breve{\tau}$ is the reverse cosmological time of $(N,g)$.

        By Proposition \ref{the duality between convex surfaces}, the dual surface $S^{*}_K$ of $S_K$ is obtained by pushing $S_K$ along orthogonal geodesics  in the future direction for a distance $\pi/2$. Then the length of a timelike curve joining $S_K$ to $S^{*}_{K}$ is at most $\pi/2$. Hence,
        \begin{equation*}
        \lim_{K\rightarrow -1}\sup\limits_{x\in S_{K}}\breve{\tau}(x)\leq\pi/2.
        \end{equation*}

        Note that $B^{-}=\{x\in N: \breve{\tau}(x)>\pi/2\}$ and $S_K$ is contained in $B^{-}$ for all $K<-1$. Therefore,
        \begin{equation}\label{eq:cosmological time of S_K}
            \lim_{K\rightarrow -1}\sup\limits_{x\in S_{K}}\breve{\tau}(x)=\pi/2.
        \end{equation}

        For any point $x\in B^{-}\setminus V_{1}$, we have $\breve{\tau}(x)>\pi/2$. By \eqref{eq:cosmological time of S_K}, there exists a future-convex surface $S_{K_x}$ of constant curvature $K_{x}\in[K_1,-1)$ such that $x$ is in the past of $S_{K_x}$. Therefore, $x\in V_{2}$. This completes the proof.
    \end{proof}

    \begin{proposition}\label{foliation of the past component}
        Let $(N,g)$ be a convex GHM AdS manifold with particles. Then $B^{-}$ admits a unique foliation by future-convex spacelike surfaces of constant curvature. Moreover, the curvature varies from $-1$ near the upper boundary component of $C(N)$ to $-\infty$ near the past singularity of $N$.
    \end{proposition}

    \begin{proof}
        By Lemma \ref{curvatures and positions of surfaces}, the future-convex spacelike surface $S_{K}$ of constant curvature $K$ (which is orthogonal to the singular lines) is unique. Moreover, $S_{K}$ and $S_{K'}$ are disjoint for all $K\not=K'\in(-\infty,-1)$. Combining this and Proposition \ref{the filling is full}, we obtain the existence and uniqueness.
    \end{proof}

        By Proposition \ref{the duality between convex surfaces}, the corresponding result of Proposition \ref{foliation of the past component} also holds for $B^{+}$ as follows.

    \begin{corollary}\label{foliation of the future component}
         Let $(N,g)$ be a convex GHM AdS manifold with particles.
         Then $B^{+}$ admits a unique foliation of past-convex spacelike surfaces of constant curvature. Moreover, the curvature varies from $-1$ near the upper  boundary component of $C(N)$ to $-\infty$ near the future singularity of $N$.
    \end{corollary}

    \emph{Proof of Theorem \ref{foliation of the GHMC AdS manifold}.} It follows directly from Proposition
    \ref{foliation of the past component} and
    Corollary \ref{foliation of the future component}.

    \begin{remark}
        It follows from Statement (2) of Proposition \ref{the convergence of the sequence of convex surfaces} and Proposition \ref{the duality between convex surfaces} that the future-convex (resp. past-convex) spacelike surface $S_{K}$ of constant curvature $K$ (which is orthogonal to the singular lines) depends continuously on $K\in(-\infty, -1)$. This implies that the (unique) foliation of $B^{-}$ (resp. $B^{+}$) is a continuous foliation.
    \end{remark}

    \section{Applications.}

    In this section, we use the results obtained above on K-surfaces in convex GHM AdS spacetimes with particles to extend to hyperbolic surfaces with cone singularities (of fixed angles less than $\pi$) a number of results concerning the landslide flow (see e.g. \cite{BMS1}). Hence we give a partial answer to the last question posed in Section 9 of \cite{BMS1}.

    Using Theorem \ref{Toulisse1} and Proposition \ref{Toulisse2}, we extend the definition of a landslide action of $S^1$ on $\mathcal{T}_{\Sigma,\theta}\times\mathcal{T}_{\Sigma,\theta}$. Moreover, as an application of Theorem \ref{foliation of the GHMC AdS manifold}, we extend to hyperbolic surfaces with cone singularities an analog of Thurston's Earthquake Theorem for the landslide flow on $\mathcal{T}_{\Sigma,\theta}\times\mathcal{T}_{\Sigma,\theta}$. Finally, we show that the relation between the AdS geometry and landslides provides more details about the parametrization map $\phi_K$.

    \subsection{The landslide action of $S^1$ on $\mathcal{T}_{\Sigma,\theta}\times\mathcal{T}_{\Sigma,\theta}$.}
        First we define the landslide transformation on $\mathcal{T}_{\Sigma,\theta}\times\mathcal{T}_{\Sigma,\theta}$.

        Let $(h,h')\in\mathfrak{M}^{\theta}_{-1}\times\mathfrak{M}^{\theta}_{-1}$, let $b:T\Sigma\rightarrow T\Sigma$ be the bundle morphism associated to $h$ and $h'$ by Corollary \ref{hyperbolic metrics and bundle morphism}, and let $\alpha\in\mathbb{R}$. Set
        \begin{equation*}
            \beta_{\alpha}=\cos{\left(\frac{\alpha}{2}\right)}E+\sin\left(\frac{\alpha}{2}\right)Jb,
        \end{equation*}
        where $E:T\Sigma\rightarrow T\Sigma$ is the identity morphism and $J$ is the complex structure induced by $h$.

        Let $h_\alpha=h(\beta_{\alpha}\bullet,\beta_{\alpha}\bullet)$ and define
        \begin{equation*}
            L_{e^{i\alpha}}(h,h'):=(h_{\alpha},h_{\alpha+\pi}).
        \end{equation*}
        In particular, we have $L_{1}(h,h')=(h,h')$, and $L_{-1}(h,h')=(h',h)$. Denote by $L^1_{e^{i\alpha}}$ (resp. $L^2_{e^{i\alpha}}$) the composition of $L_{e^{i\alpha}}$ with the projection on the first (resp. second) factor.

        \begin{proposition}
            For all $\alpha\in\mathbb{R}$, $h_{\alpha}$ is a hyperbolic metric with cone singularities of the same angles as $h$.
        \end{proposition}

        \begin{proof}
            It can be checked (as in the proof of Lemma 3.2 in \cite{BMS1}) that $d^{\nabla}\beta_{\alpha}=0$  and $\det(\beta_{\alpha})=1$, where $\nabla$ is the Levi-Civita connection of $h$. By Theorem \ref{computation of connection and curvature}, the Levi-Civita connection $\nabla^{\alpha}$ of $h_{\alpha}$ is given by $\nabla^{\alpha}_{u}v=\beta_{-\alpha}\nabla_{u}(\beta_{\alpha}v)$. The curvature of $h_{\alpha}$ outside the singular locus is
            \begin{equation}\label{hyperbolic metric}
                K_{\alpha}=\frac{K_{h}}{\det(\beta_\alpha)}=-1.
            \end{equation}

            Note that $\beta_{\alpha}=\cos{(\frac{\alpha}{2})}E+\sin(\frac{\alpha}{2})Jb$, where $b:T\Sigma\rightarrow T\Sigma$ is a bundle morphism associated to $(h,h')$ by Corollary \ref{hyperbolic metrics and bundle morphism}.
            In particular, both eigenvalues of $b$ tend to 1 at the marked points of $\Sigma$.

             We can choose a suitable orthonormal frame $(e_1,e_2)$ such that $h(e_{i},e_{j})=\delta_{ij}$ for $i,j=1,2$, and $b(e_{1})=k_1e_1$, $b(e_{2})=k_2e_2$, where $k_1,k_2>0$. In this frame, we have the following expressions:
            \begin{equation*}
                J=\left(
                        \begin{array}{cc}
                            0 & -1 \\
                            1 & 0\\
                        \end{array}
                  \right),\qquad
                b=\left(
                        \begin{array}{cc}
                            k_1 & 0 \\
                            0 & k_2\\
                        \end{array}
                 \right),\qquad
                \beta_\alpha=\left(
                                \begin{array}{cc}
                                    \cos(\frac{\alpha}{2}) & -k_2\sin(\frac{\alpha}{2}) \\
                                    k_1\sin(\frac{\alpha}{2}) & \cos(\frac{\alpha}{2})\\
                                \end{array}
                            \right).
            \end{equation*}

            A direct computation shows that the matrix of $h_{\alpha}=h(\beta_{\alpha}\bullet,\beta_{\alpha}\bullet)$ in this frame is
            \begin{equation*}
                h_\alpha=\left(
                            \begin{array}{cc}
                                k_1^2\sin^{2}(\frac{\alpha}{2})+\cos^{2}(\frac{\alpha}{2})& \frac{1}{2}(k_1-k_2)\sin(\frac{\alpha}{2}) \\\\
                                \frac{1}{2}(k_1-k_2)\sin(\frac{\alpha}{2})& \sin^{2}(\frac{\alpha}{2})+k_2^2\cos^{2}(\frac{\alpha}{2}) \\
                            \end{array}
                        \right).
         \end{equation*}

            Note that $k_1,k_2$ tend to 1 at the marked points of $\Sigma$. Therefore, $h_{\alpha}$ tends to $h$ at the marked points. Combining this with \eqref{hyperbolic metric}, we obtain that $h_{\alpha}$ is a hyperbolic metric with cone singularities, with the same cone angles as $h$.
        \end{proof}

        The following lemma shows that the landslide flow is well-defined on $\mathcal{T}_{\Sigma,\theta}\times\mathcal{T}_{\Sigma,\theta}$.

        \begin{lemma}
            Let $(\tau,\tau')\in\mathcal{T}_{\Sigma,\theta}\times\mathcal{T}_{\Sigma,\theta}$, and let $(h,h')$, $(\bar{h},\bar{h}')$ be two normalized representatives of $(\tau,\tau')$. Then $h_{\alpha}=L^{1}_{e^{i\alpha}}(h,h')$ and  $\bar{h}_{\alpha}=L^{1}_{e^{i\alpha}}(\bar{h},\bar{h'})$ are isotopic in $\mathcal{T}_{\Sigma,\theta}$ for all $\alpha\in\mathbb{R}$.
        \end{lemma}

        \begin{proof}
            By definition, for all $\alpha\in\mathbb{R}$,
            \begin{equation*}
                h_\alpha=h(\beta_\alpha\bullet,\beta_\alpha\bullet),\qquad
                \bar{h}_\alpha=\bar{h}(\bar{\beta}_\alpha\bullet,\bar{\beta}_\alpha\bullet),
            \end{equation*}
            where $\beta_{\alpha}=\cos(\frac{\alpha}{2})+\sin(\frac{\alpha}{2})Jb$, and $\bar{\beta}_{\alpha}=\cos(\frac{\alpha}{2})+\sin(\frac{\alpha}{2})\bar{J}\,\bar{b}$. Here $b$ (resp. $\bar{b}$) is the bundle morphism associated to $h,h'$ (resp. $\bar{h},\bar{h'}$) by Corollary \ref{hyperbolic metrics and bundle morphism}. $J$ (resp. $\bar{J}$) is the complex structure induced by $h$ (resp. $\bar{h}$).

            By Remark \ref{rk: normalized representatives}, there exists a diffeomorphism $f$ from $\Sigma$ to $\Sigma$, which is isotopic to the identity (the isotopy fixes each marked point) such that $\bar{h}=f^{*}h$, $\bar{h'}=f^{*}h'$. Using a similar argument as in the proof of Lemma \ref{well-define of the map}, we can prove that
            \begin{equation}\label{eq:b and J}
                \bar{b}=(df)^{-1}b(df)~,\qquad
                \bar{J}=(df)^{-1}J(df)~,
            \end{equation}
            where $df$ is the differential of $f$. Applying \eqref{eq:b and J} to the expression of $\bar{\beta_{\alpha}}$, we obtain that
            \begin{equation}\label{eq:beta_varphi}
                \bar{\beta}_{\alpha}
                =(df)^{-1}(\cos(\frac{\alpha}{2})+\sin(\frac{\alpha}{2})Jb)(df)
                =(df)^{-1}\beta_\alpha (df).
            \end{equation}
            Substituting \eqref{eq:beta_varphi} and $\bar{h}=f^{*}h$ into $\bar{h}_\alpha=\bar{h}(\bar{\beta}_\alpha\bullet,\bar{\beta}_\alpha\bullet)$, we see that $\bar{h}_\alpha=f^{*}(h_{\alpha})$. This implies that $\bar{h}_{\alpha}$ is isotopic to $h_{\alpha}$ for all $\alpha\in\mathbb{R}$.
        \end{proof}

        \begin{remark}
            For simplicity, henceforth we denote by $(h,h')$ both a pair of normalized metrics in $\mathfrak{M}^{\theta}_{-1}\times\mathfrak{M}^{\theta}_{-1}$ and its equivalence class in $\mathcal{T}_{\Sigma,\theta}\times\mathcal{T}_{\Sigma,\theta}$.
        \end{remark}

        \begin{definition}\label{defintion of the landslide}
            For all $\alpha\in\mathbb{R}$, the map $L_{e^{i\alpha}}:\mathcal{T}_{\Sigma,\theta}\times\mathcal{T}_{\Sigma,\theta} \rightarrow\mathcal{T}_{\Sigma,\theta}\times\mathcal{T}_{\Sigma,\theta}$ sending an element $(h,h')\in\mathcal{T}_{\Sigma,\theta}\times\mathcal{T}_{\Sigma,\theta}$ to the pair of isotopy classes of $h_{\alpha}$, $h_{\alpha+\pi}$ is well-defined. We call $L_{e^{i\alpha}}$ the landslide (transformation) of parameter $\alpha$.
        \end{definition}

        Note that the argument for the flow property of landslides on the product of two copies of the Teichm\"uller space of a closed surface (see Theorem 1.8 in \cite{BMS1}) can be directly applied to the case with cone singularities. It leads to the following proposition.

        \begin{proposition}
            The landslide $L_{e^{i\alpha}}$ given by Definition \ref{defintion of the landslide} is a flow: for any $\alpha,\alpha'\in\mathbb{R}$,
            \begin{equation*}
                L_{e^{i\alpha'}}\circ L_{e^{i\alpha}}=L_{e^{i(\alpha+\alpha')}}.
            \end{equation*}
            In other words, the map $L:\mathcal{T}_{\Sigma,\theta}\times\mathcal{T}_{\Sigma,\theta}\times S^{1}\rightarrow\mathcal{T}_{\Sigma,\theta}\times\mathcal{T}_{\Sigma,\theta}$ associating to $(h,h',e^{i\alpha})$ the image $L_{e^{i\alpha}}(h,h')$ defines an action of $S^{1}$ on $\mathcal{T}_{\Sigma,\theta}\times\mathcal{T}_{\Sigma,\theta}$. We call $L$ the landslide flow, or the landslide action on $\mathcal{T}_{\Sigma,\theta}\times\mathcal{T}_{\Sigma,\theta}$.
        \end{proposition}

    \subsection{The extension of Thurston's Earthquake Theorem.}

        In this section we extend to hyperbolic surfaces with cone singularities (of fixed angles less than $\pi$) an analog of the Earthquake Theorem, already proved for the landslide flow on non-singular hyperbolic surfaces in \cite{BMS1}. To prove this theorem, we give the following lemma, as a generalization of Lemma 1.9 in \cite{BMS1} to the case with cone singularities.

\begin{lemma}\label{Lemma1 for extension of Thurston's earthquake theorem}
            Let  $(h,h')\in\mathfrak{M}^{\theta}_{-1}\times\mathfrak{M}^{\theta}_{-1}$ be a pair of normalized metrics and let $\alpha\in(0,\pi)$. Then there exists a unique GHM AdS spacetime $(N,g)$ with particles which contains a future-convex spacelike surface orthogonal to the singular lines with the induced metric $I_{\alpha}=\cos^{2}(\frac{\alpha}{2})h$ and the third fundamental form $III_{\alpha}=\sin^{2}(\frac{\alpha}{2})h'$. Moreover, $L^{1}_{e^{i\alpha}}(h,h')$ and $L^{1}_{e^{-i\alpha}}(h,h')$ are the left and right metrics of $(N,g)$, respectively.
\end{lemma}

        \begin{proof}
            Note that $\cos^2(\frac{\alpha}{2}),\sin^2(\frac{\alpha}{2})\in(0,1)$.  The first statement is a direct consequence of Lemma \ref{definition of the map} applied with $K=-1/\cos^2(\frac{\alpha}{2})$ and $K^{*}=-1/\sin^2(\frac{\alpha}{2})$.

            Denote by $B_{\alpha}$ the shape operator of the future-convex spacelike surface of constant curvature $K$ in $(N,g)$ and denote by $J_{\alpha}$ the complex structure of $I_{\alpha}$. A simple computation shows that $B_{\alpha}=\tan(\frac{\alpha}{2})b$, where $b$ is the bundle morphism associated to $h,h'$ by Corollary \ref{hyperbolic metrics and bundle morphism}, and $J_{\alpha}=J$, where $J$ is the complex structure of $h$. By the extension of Mess' parametrization (see Theorem 1.4 in \cite{BS}), the left and right metrics of $(N,g)$ can be expressed as
            \begin{equation}\label{left metric}
                \mu_{l}
                =I_{\alpha}((E+J_{\alpha}B_{\alpha})\bullet,(E+J_{\alpha}B_{\alpha})\bullet),\qquad
                \mu_{r}
                =I_{\alpha}((E-J_{\alpha}B_{\alpha})\bullet,(E-J_{\alpha}B_{\alpha})\bullet).
            \end{equation}

            Substituting $B_{\alpha}=\tan(\frac{\alpha}{2})b$ and $J_{\alpha}=J$ into \eqref{left metric}, we obtain that
            \begin{equation*}
                \mu_{l}
                =h(\beta_{\alpha}\bullet,\beta_{\alpha}\bullet)
                =L^{1}_{e^{i\alpha}}(h,h').
            \end{equation*}

            Similarly, we can prove that the right metric of $(N,g)$ is
            \begin{equation*}
                \mu_{r}
                =h(\beta_{-\alpha}\bullet,\beta_{-\alpha}\bullet)
                =L^{1}_{e^{-i\alpha}}(h,h').
            \end{equation*}
            This completes the proof.
        \end{proof}

        \begin{corollary}\label{corolalry about landslide}
            Let $(\mu_{l},\mu_{r})\in\mathcal{T}_{\Sigma,\theta}\times\mathcal{T}_{\Sigma,\theta}$ and $\alpha\in(0,\pi)$. There exists a unique $(h,h')\in\mathcal{T}_{\Sigma,\theta}\times\mathcal{T}_{\Sigma,\theta}$ such that
            \begin{equation*}
                L^{1}_{e^{i\alpha}}(h,h')=\mu_{l},\qquad
                L^{1}_{e^{-i\alpha}}(h,h')=\mu_{r}.
            \end{equation*}
        \end{corollary}

        \begin{proof}
            Given $\mu_{l}$ and $\mu_{r}$, by the extension of Mess' Parametrization (see Theorem 1.4 in \cite{BS}), there exists a unique convex GHM AdS manifold $(N,g)$ with particles which has the left and right metrics $\mu_{l}$ and $\mu_{r}$. By Theorem \ref{foliation of the GHMC AdS manifold}, $(N,g)$ contains a unique future-convex surface $S_K$ of constant curvature $K=-1/\cos^{2}(\frac{\alpha}{2})$. Denote by $I$ and $III$ the first and third fundamental form on $S_K$. Then $III$ has constant curvature $K^{*}=-1/\sin^{2}(\frac{\alpha}{2})$. Set $h=|K|I$ and $h'=|K^{*}|III$.
            It can be checked that $(h,h')$ is a pair of normalized metrics. It follows from Lemma \ref{Lemma1 for extension of Thurston's earthquake theorem} that $L^{1}_{e^{i\alpha}}(h,h')=\mu_{l}$, $L^{1}_{e^{-i\alpha}}(h,h')=\mu_{r}$. This shows the existence.

            Now we show the uniqueness. Suppose $(\bar{h},\bar{h'})\in\mathcal{T}_{\Sigma,\theta}\times\mathcal{T}_{\Sigma,\theta}$ is another pair such that
            \begin{equation}\label{eq:uniqueness}
                L^{1}_{e^{i\alpha}}(\bar{h},\bar{h'})=\mu_{l},\qquad
                L^{1}_{e^{-i\alpha}}(\bar{h},\bar{h'})=\mu_{r}.
            \end{equation}

            By Lemma \ref{Lemma1 for extension of Thurston's earthquake theorem}, there exists a unique GHM AdS manifold $(\bar{N},\bar{g})$ with particles which contains a future-convex spacelike surface orthogonal to the singular lines, with the induced metric $\cos^{2}(\frac{\alpha}{2})\bar{h}$ and the third fundamental form $\sin^{2}(\frac{\alpha}{2})\bar{h'}$. Moreover, by \eqref{eq:uniqueness}, the left and right metrics of $(\bar{N},\bar{g})$ are $\mu_{l}$ and $\mu_{r}$, respectively. The extension of Mess' parametrization implies that $(\bar{N},\bar{g})$ is $(N,g)$ up to isotopy. The uniqueness in Theorem \ref{foliation of the GHMC AdS manifold} shows that $(\bar{h},\bar{h'})=(h,h')$ in $\mathcal{T}_{\Sigma,\theta}\times\mathcal{T}_{\Sigma,\theta}$.
        \end{proof}

        Now we are ready to prove the extension of Thurston's Earthquake Theorem to the case with cone singularities, which generalizes Theorem 1.14 in \cite{BMS1}.

        \begin{theorem} \label{tm:landslides}
            Let $h_{1},h_{2}\in\mathcal{T}_{\Sigma,\theta}$ and let $e^{i\alpha}\in S^{1}\setminus\{1\}$. Then there exists a unique $h_{1}'\in\mathcal{T}_{\Sigma,\theta}$ such that $L^{1}_{e^{i\alpha}}(h_{1},h_{1}')=h_{2}$.
        \end{theorem}

        \begin{proof}
            First we show the existence. Corollary \ref{corolalry about landslide} applied with $\mu_l=h_{2}$, $\mu_r=h_{1}$ and $\varphi=\alpha/2$ shows that there exists a unique $(h_{0},h_{0}')\in\mathcal{T}_{\Sigma,\theta}\times\mathcal{T}_{\Sigma,\theta}$ such that $L^{1}_{e^{i\varphi}}(h_0,h_0')=h_{2}$ and $L^{1}_{e^{-i\varphi}}(h_0,h_0')=h_{1}$. Set $h_{1}'=L^{2}_{e^{-i\varphi}}(h_0,h_0')$. Then we get
            \begin{equation*}
                L^{1}_{e^{i\alpha}}(h_{1},h_{1}')
                =L^{1}_{e^{i2\varphi}}(L_{e^{-i\varphi}}(h_0,h_0'))
                =h_{2}.
            \end{equation*}

            Assume that $\bar{h}'_1$ is another element in $\mathcal{T}_{\Sigma,\theta}$ such that $L^{1}_{e^{i\alpha}}(h_{1}$,$\bar{h}'_1)=h_{2}$. Set $(h,h')=L_{e^{i\alpha/2}}(h_1,\bar{h}'_1)$. By computation, we have $L^{1}_{e^{i\alpha/2}}(h,h')=h_2$, and $L^{1}_{e^{-i\alpha/2}}(h,h')=h_1$. The uniqueness in Corollary \ref{corolalry about landslide} implies that $(h_{0},h_{0}')=(h,\bar{h}')$. Hence $\bar{h}'_1=L^{2}_{e^{-i\varphi}}(h,h')=L^{2}_{e^{-i\varphi}}(h_0,h_0')=h_1'$. This completes the proof.
        \end{proof}

    \subsection{The landslide flow in terms of harmonic maps.}

        Recall that in the non-singular case, landslides can also be defined in terms of multiplication of the Hopf differential of harmonic maps by complex numbers of modulus $1$ (see Definition 1.2 in \cite{BMS1}).

        Consider a map $\Phi:\mathcal{T}_{\Sigma,\theta}\to \mathcal{QD}_{c}(\Sigma)$, which associates to $g\in\mathcal{T}_{\Sigma,\theta}$ the Hopf differential (with respect to the conformal structure $c$) of the harmonic map $u_{c,\bar{g}}$ from $(\Sigma,c)$ to $(\Sigma,\bar{g})$ isotopic to the identity, where $\bar{g}$ is a representative of $g$. By the uniqueness in Theorem \ref{existence and uniqueness of harmonic maps} and the fact that harmonic maps remain harmonic when composed from the left with isometries, $\Phi$ is well-defined (i.e. independent of the choice of the representatives of $g$).

        For simplicity, we use the same notation for both $g\in\mathcal{T}_{\Sigma,\theta}$ and its representative in the proof of the following proposition. The argument is similar to that for Theorem 3.1 in \cite{Wolf}.

\begin{proposition}
            \label{homeomorphic between Teichmuller space and quadratic differential}
The map $\Phi$ is a homeomorphism.
\end{proposition}

\begin{proof}
            Observe that $\mathcal{T}_{\Sigma,\theta}$ and $\mathcal{QD}_{c}(\Sigma)$ are both $6g-6+2n$-dimensional cells. By Brouwer's Invariance of Domain Theorem, it suffices to show that $\Phi$ is continuous, injective and proper.

            The continuity is obvious, since the harmonic maps $u_{c,g}$ vary smoothly with respect to the target metric $g$ (see Theorem \ref{existence and uniqueness of harmonic maps}).

            For the injectivity of $\Phi$, we use the maximum principle as applied in \cite[Theorem 3.1]{Wolf}. Suppose that $g_1,g_2\in\mathcal{T}_{\Sigma,\theta}$ satisfy that $\Phi(g_1)=\Phi(g_2)$. Denote $\Phi(g_i)=\Phi_i$ ($i$ always takes values in $\{1,2\}$ in this proof), so that $\Phi_1=\Phi_2$. Let $z$ be a conformal coordinate on $(\Sigma,c)$. Set $c=g_0=\sigma(z)|dz|^2$, $g_{i}=\rho_i(u_i(z))|du_i|^2$, where $u_i=u_{c,g_i}$. By computation, we obtain that $\Phi_i=\rho_i(u_i(z))(u_{i})_{z}\overline{(u_{i})_{\bar{z}}}$. Set $H_i=\sigma^{-1}(z)\rho_{i}(u_i(z))|(u_{i})_z|^2$, and $L_i=\sigma^{-1}(z)\rho_{i}(u_i(z))|(u_{i})_{\bar{z}}|^2$. We have the following quantities (see \cite[Section 2]{Wolf}):
            \begin{enumerate}[(a)]
                \item The energy density $e_i=H_i+L_i$.
                \item The Jacobian $J_i=H_i-L_i>0$.
                \item The norm of the quadratic differential $|\Phi_i|^2/{\sigma}^2=H_iL_i$.
                \item The Beltrami differential $\nu_i=(u_i)_{\bar{z}}/(u_i)_{z}=\overline{\Phi_i}/\sigma H_i$.
                \item The pull-bak metric of $\rho_i$ by $u_i$ is ${u_i}^{*}(g_i)=2\Re(\Phi_idz^2)+\sigma e_i dzd\bar{z}$.
           \end{enumerate}

           Set $h_i=\log H_i$ and $\Delta=4{\sigma}^{-1}\partial^{2}_{z\bar{z}}$. It is well-known that the following identity (see \cite{SY}) holds:
           \begin{equation*}
                \Delta h_i=2(H_i-L_i-1)=2(e^{h_i}-\sigma^{-2}|\Phi_i|^2e^{-h_i}-1).
           \end{equation*}

           We claim that $h_1=h_2$. Indeed, if there exists a point at which $h_1>h_2$, there exists a maximum point $x_0\in\Sigma$ of $h_1-h_2$, at which $h_1-h_2>0$. To see that this cannot be, note that there exists a neighbourhood $U$ of the cone singularity $p_k$ of angle $\beta_k$, such that $\sigma(z)=e^{\lambda}|z|^{2(\beta_k-1)}$, $\rho_i(u_i(z))=e^{\zeta_i}|u_i(z)|^{2(\beta_k-1)}$, where $\beta_k=\theta_k/(2\pi)$
           $z$ is the conformal coordinate centered at $p_k$, and $\lambda$, $\zeta_i$  are continuous functions on $U$. Moreover, $u_i(z)=\xi_{i}z + r^{1+\varepsilon}f_i(z)$, where $\xi_i\in\mathbb{C}\setminus\{0\}$, $r=|z|$, $\varepsilon>0$ and $f_{i}\in\chi^{2,\gamma}_{b}(U)\cap\mathcal{C}^{2}(U)$ (see \cite[Section 4.2]{Toulisse2}). It is also computed in \cite[Section 4.2]{Toulisse2} that
           \begin{equation}\label{eq:H and L}
           \begin{split}
           H_i&=\sigma^{-1}(z)e^{\zeta_i}|\xi_{i}|^{2\beta_k}r^{2(\beta_k-1)}
           (1+O(r^{\varepsilon})),\\
           L_i&=\sigma^{-1}(z)e^{\zeta_i}|\xi_{i}|^{2(\beta_k-1)}r^{2(\beta_k-1)+2\varepsilon}
           |\bar{L}(f_{i})|^{2}(1+O(r^{\varepsilon})),
           \end{split}
           \end{equation}
           where the operator $\bar{L}=\frac{r}{2\bar{z}}((1+\varepsilon)Id+r\partial_r+i\partial_{\alpha})$ and $z=re^{i\alpha}$.

           Substituting $\sigma(z)=e^{\lambda}|z|^{2(\beta_k-1)}$ into \eqref{eq:H and L}, we obtain that
           \begin{equation*}
           \begin{split}
           H_i&=e^{\zeta_i-\lambda}|\xi_{i}|^{2\beta_k}(1+O(r^{\varepsilon})),\\
           L_i&=e^{\zeta_i-\lambda}|\xi_{i}|^{2(\beta_k-1)}r^{2\varepsilon}
           |\bar{L}(f_{i})|^{2}(1+O(r^{\varepsilon})).
           \end{split}
           \end{equation*}
           Note that we can 
           make a completion of the punctured surfaces $(\Sigma_{\mathfrak{p}}, g_0)$ and ($\Sigma_{\mathfrak{p}}, g_i)$ by directly adding the set $\mathfrak{p}$.
           Hence $h_1-h_2=\log(H_1/H_2)$ can be viewed as a $\mathcal{C}^{2}$ function on a compact surface $\Sigma$ and has a maximum point. At this point, we have
           \begin{equation*}
                0\geq\Delta(h_1-h_2)=2\{(e^{h_1}-e^{h_2})-\sigma^{-2}|\Phi|^2(e^{-h_1}-e^{-h_2})\}>0.
           \end{equation*}

           This implies that $h_1\leq h_2$. Symmetrically, $h_2\leq h_1$. Hence, $h_1=h_2$, which implies that $H_1=H_2$. By equality (c) and (a), $L_1=L_2$ and $e_1=e_2$. Combined with equality (e), we get ${u_1}^{*}(g_1)={u_2}^{*}(g_2)$. Note that $u_1, u_2$ are isotopic to identity, then $g_1=g_2\in\mathcal{T}_{\Sigma,\theta}$.

            Define the map $E:\mathcal{T}_{\Sigma,\theta}\to \mathbb{R}$ as $E(g)=E(u_{c,g})$. To show the properness of $\Phi$, we first state the fact that $||\Phi(g)||\to\infty$ iff $E(g)\to\infty$. Indeed, applying equalities (b),(c),(d) and the Gauss-Bonnet formula for surfaces with cone singularities (see e.g. \cite[Proposition 1]{Troyanov}):
            \begin{equation*}
                \int J\sigma dzd\bar{z}
                =\area_g(\Sigma)
                =-2\pi\{\chi(\Sigma)+\sum\limits_{k=1}^{n_0}(\theta_k/{2\pi}-1)\}
                :=-2\pi\chi(\Sigma,\theta),
            \end{equation*}
            as in \cite[Theorem 3.1]{Wolf}, we have
            \begin{equation*}
                    \int H\sigma dzd\bar{z}+2\pi\chi(\Sigma,\theta)
                \leq\int L\sigma dzd\bar{z}
                \leq\int |\Phi(g)|dzd\bar{z}
                \leq\int H\sigma dzd\bar{z}
                \leq\int L\sigma dzd\bar{z}-2\pi\chi(\Sigma,\theta).
            \end{equation*}

            Adding the first two and last two integrals and applying equality (a), we obtain
            \begin{equation*}
                    E(g)+2\pi\chi(\Sigma,\theta)
                \leq 2\int|\Phi(g)|dzd\bar{z}
                \leq E(g)-2\pi\chi(\Sigma,\theta).
            \end{equation*}

            Now we are left to show that $E$ is proper. That is, to show that $B=\{g\in\mathcal{T}_{\Sigma,\theta}:E(g)<C_0\}$ is compact. By Lemma \ref{the compactness of some subset of Teichmuller space}, it suffices to show that
            \begin{equation}
                \label{ineq for property}
                    \ell_{\gamma}(g)\leq C\ell_{\gamma}(g_0),
            \end{equation}
            for all $g\in B$ and all simple closed curves $\gamma$ on $\Sigma$.

            By Lemma \ref{Collar lemma for hyperbolic cone-surfaces}, there exists a uniform lower bound for the injectivity radius of the singularities over $\mathcal{T}_{\Sigma,\theta}$. Denote by $\overline{\Sigma}^{g_0}$ (resp. $\overline{\Sigma}^{g}$) 
            the completion of $(\Sigma_{\mathfrak{p}},g_0)$ (resp. $(\Sigma_{\mathfrak{p}},g)$) by adding the set $\mathfrak{p}$ and denote by $\inj(g_0)$ the injectivity radius of $\overline{\Sigma}^{g_0}$. Then $\inj(g_0)>0$. Let $c_1(g_0)=\min\{1,(\inj(g_0))^2\}$. The Courant-Lebesgue Lemma (see \cite[Proposition 3]{Wolf} and \cite[Lemma 3.1]{Jost}) can be applied to the harmonic map $u:\overline{\Sigma}^{g_0}\to \overline{\Sigma}^{g}$, that is, for any $x_1, x_2\in\overline{\Sigma}^{g_0}$ with $d_{g_0}(x_1,x_2)<\delta<c_1(g_0)$,
            \begin{equation*}
                d_{g}(u(x_1),u(x_2))<4\sqrt{2}\pi{C_{0}}^{1/2}(\log(1/\delta))^{-1/2}.
            \end{equation*}
            This implies \eqref{ineq for property}, where $C$ depends on $g_0$ and $C_0$. The proof is complete.
\end{proof}

        Proposition \ref{homeomorphic between Teichmuller space and quadratic differential} shows that given a meromorphic quadratic differential $q\in\mathcal{QD}_{c}(\Sigma)$ with at most simple poles at singularities, there exists a unique $h\in\mathcal{T}_{\Sigma,\theta}$ such that the identity map $id:(\Sigma,c)\rightarrow(S,h)$ is harmonic with Hopf differential $q$.

        This statement, combined with Lemma \ref{existence and uniqueness of harmonic maps}, makes it possible to generalize the definition of the landslide flow in terms of harmonic maps to hyperbolic surfaces with cone singularities as follows.

        \begin{definition}\label{definition of R}
            Let $c,h\in\mathcal{T}_{\Sigma,\theta}$ and let $e^{i\alpha}\in S^1$. Define $R_{c,\alpha}(h)$ as the (unique) metric $h^{\alpha}\in\mathcal{T}_{\Sigma,\theta}$ such that if $f:(\Sigma,c)\rightarrow(\Sigma,h)$ and $f^{\alpha}:(\Sigma,c)\rightarrow(\Sigma,h^{\alpha})$ are the harmonic maps isotopic to the identity (fixing each marked point), then $\Phi(f^{\alpha})=e^{i\alpha}\Phi(f)$.
        \end{definition}

        Let $h, h'\in\mathcal{T}_{\Sigma,\theta}$. Recall that if $h_\alpha$ is used to denote $L^1_{e^{i\alpha}}(h,h')$, then $L_{e^{i\alpha}}(h,h')=(h_{\alpha},h_{\alpha+\pi})$. Denote by $c_{\alpha}$ the conformal structure of the metric $h_{\alpha}+m^{*}_{\alpha}(h_{\alpha+\pi})$, where $m_{\alpha}:(\Sigma,h_{\alpha})\rightarrow(\Sigma,h_{\alpha+\pi})$ is the unique minimal Lagrangian map isotopic to the identity, which is called the \emph{center} of $(h_{\alpha},h_{\alpha+\pi})$. Applying the analogous argument as Theorem 1.10 in \cite{BMS1} to the case with cone singularities, we have the following proposition.

        \begin{proposition} \label{center proposition}
            Let $h, h'\in\mathcal{T}_{\Sigma,\theta}$ and let $c_{\alpha}$ be the center of $(h_{\alpha},h_{\alpha+\pi})$. Then
            \begin{enumerate}[(1)]
                \item  The identity $id:(\Sigma,h_{\alpha})\rightarrow(\Sigma,h_{\alpha+\pi})$ is minimal Lagrangian.
                \item  $c_{\alpha}$ is independent of $\alpha$ --- we denote it by $c$.
                \item  For any $\alpha\in\mathbb{R}$, $\Phi(f_{\alpha})=e^{i\alpha}\Phi(f)$, where $f_{\alpha}:(\Sigma,c)\rightarrow(\Sigma,h_{\alpha})$ is the unique harmonic map isotopic to the identity..
            \end{enumerate}
        \end{proposition}

        The following corollary is a direct consequence of Definition \ref{definition of R} and Proposition \ref{center proposition}.

        \begin{corollary}\label{equivalent definiton}
            Let $(h,h')\in\mathcal{T}_{\Sigma,\theta}\times\mathcal{T}_{\Sigma,\theta}$ be a normalized representative, and let $c$ be the conformal class of $h+h'$. Then for any $e^{i\alpha}\in S^1$, we have
            \begin{equation*}
                L_{e^{i\alpha}}(h,h')=(R_{c,\alpha}(h),R_{c,\alpha+\pi}(h)).
            \end{equation*}
        \end{corollary}

    \subsection{An application of the landslide flow.}

        We now go in the reverse direction, and use the properties of the landslide flow to obtain new results on the geometry of $K$-surfaces in convex GHM AdS spacetimes with particles. We first state a lemma on landslides on hyperbolic surfaces with cone singularities, and then use it to obtain Theorem \ref{tm:K1K2} below on $K$-surfaces.

        \begin{lemma}\label{lemma for the injectivity of the landslide}
            Let $(h,h')\in\mathcal{T}_{\Sigma,\theta}\times\mathcal{T}_{\Sigma,\theta}$ be a normalized representative. Define the map $L_{\bullet}(h,h'):S^{1}\rightarrow\mathcal{T}_{\Sigma,\theta}\times\mathcal{T}_{\Sigma,\theta}$ by associating $L_{e^{i\alpha}}(h,h')$ to $e^{i\alpha}\in S^1$. Then the following two statements hold:
            \begin{enumerate}[(1)]
                \item  If $h\not= h'$, then the map $e^{i\alpha}\mapsto L_{e^{i\alpha}}(h,h')$ is injective.
                \item  If $h=h'$, then this map $e^{i\alpha}\mapsto L_{e^{i\alpha}}(h,h')$ is constant, that is, $L_{e^{i\alpha}}(h,h')=(h,h)$ for all $e^{i\alpha}\in S^1$.
            \end{enumerate}
        \end{lemma}

        \begin{proof}
            First we show the first statement. Assume that $h\not=h'$ and $L_{e^{i\alpha_1}}(h,h')=L_{e^{i\alpha_2}}(h,h')$. By Corollary \ref{equivalent definiton}, we have
            \begin{equation}\label{eq for the injectivity}
                R_{c,\alpha_1}(h)=R_{c,\alpha_2}(h),\qquad
                \Phi(f^{\alpha_{i}})=e^{i\alpha_{i}}\Phi(f),
            \end{equation}
            for $i=1,2$, where $f:(\Sigma,c)\rightarrow(\Sigma,h)$ and $f^{\alpha_{i}}:(\Sigma,c)\rightarrow(\Sigma, R_{c,\alpha_i}(h))$ are the (unique) harmonic maps isotopic to the identity, $c$ is the conformal structure of $h+h'$. Moreover, $\eqref{eq for the injectivity}$ implies $\Phi(f^{\alpha_{1}})=\Phi(f^{\alpha_{2}})$, that is,
            \begin{equation*}
                e^{i(\alpha_1-\alpha_2)}\Phi(f)=0.
            \end{equation*}
            Note that $\Phi(f)\not=0$ since $h\not=h'$. This implies that $\alpha_1=\alpha_2$.

            Assume that $h=h'$, then $c$ is the conformal structure of $h$. It follows that the harmonic map $f:(\Sigma,c)\rightarrow (\Sigma,h)$ isotopic to the identity is exactly the identity by choosing the representative metric $h$ of $c$. Hence, $\Phi(f)=0$ and  $\Phi(f^{\alpha})=e^{i\alpha}\Phi(f)=0$ for all $\alpha\in S^1$. By Definition \ref{definition of R} and Corollary \ref{equivalent definiton}, we obtain
            \begin{equation*}
                L_{e^{i\alpha}}(h,h')=(R_{c,\alpha}(h),R_{c,\alpha+\pi}(h))= (h,h),
            \end{equation*}
            for all $e^{i\alpha}\in S^1$.
        \end{proof}

        \begin{theorem} \label{tm:K1K2}
            Let $(N,g)\in\mathcal{GH}_{\Sigma,\theta}$ and $K_1,K_2\in(\infty,-1)$. Then the following two statements are equivalent:
            \begin{enumerate}[(1)]
                \item  The preimages under $\phi_{K_1}$ and $\phi_{K_2}$ of $(N,g)$ are the same point in $\mathcal{T}_{\Sigma,\theta}\times\mathcal{T}_{\Sigma,\theta}$.
                \item  $K_1=K_2$ or $(N,g)$ is Fuchsian.
            \end{enumerate}
        \end{theorem}

        \begin{proof}
            First we show Statement (1) implies Statement (2). Denote by $(h,h')$ the same preimage under $\phi_{K_1}$ and $\phi_{K_2}$ of $(N,g)$. Let $\alpha_{i}\in(0,\pi)$ such that $K_{i}=-1/\cos^{2}\alpha_i$ for $i=1,2$. From the definition of $\phi_{K_i}$, $(N,g)$ contains a future-convex spacelike surface $S_{K_{i}}$ orthogonal to the singular lines, with the induced metric $(1/|K_{i})|h$ and the third fundamental form $(1/|K^{*})|h'$, where $K_{i}^{*}=-K_{i}/(1+K_{i})=-1/\sin^{2}\alpha_i$ for $i=1,2$. Apply Lemma \ref{Lemma1 for extension of Thurston's earthquake theorem} with $(h,h')\in\mathcal{T}_{\Sigma,\theta}\times\mathcal{T}_{\Sigma,\theta}$ and $\alpha_{i}\in(0,\pi)$, the left and right metrics of $(N,g)$ are respectively
            \begin{equation}\label{eq:left and right metrics}
                \mu_{l}=L^{1}_{e^{i\alpha_{1}}}(h,h')=L^{1}_{e^{i\alpha_{2}}}(h,h'),\qquad
                \mu_{r}=L^{1}_{e^{-i\alpha_{1}}}(h,h')=L^{1}_{e^{-i\alpha_{2}}}(h,h').
            \end{equation}

            We claim that if $(N,g)$ is not Fuchsian, then $h\not= h'$. Otherwise, by \eqref{eq:left and right metrics} and Statement (2) of Lemma \ref{lemma for the injectivity of the landslide}, $h=h'$ implies that $\mu_{l}=\mu_{r}$ and hence $(N,g)$ is Fuchsian. This produces contradiction. By Statement (1) of Lemma \ref{lemma for the injectivity of the landslide}, we have $\alpha_{1}=\alpha_{2}$. This implies that $K_1=K_2$.

            Now it suffices to prove that Statement (2) implies Statement (1). It follows immediately if $K_1=K_2$ from the injectivity of the map $\phi_{K_1}=\phi_{K_2}$. If $(N,g)$ is Fuchsian, denote by $(h_1,h_1')$ and $(h_2,h_2')$ the preimage under the maps $\phi_{K_1}$ and $\phi_{K_2}$ of $(N,g)$. Note that $\alpha_{i}\in(0,\pi)$. By Lemma \ref{Lemma1 for extension of Thurston's earthquake theorem}, we have
            \begin{equation*}
                \mu_{l}=L^{1}_{e^{i\alpha_{1}}}(h_1,h_1')=L^{1}_{e^{-i\alpha_{1}}}(h_1,h_1')=\mu_{r},\qquad
                \mu_{l}=L^{1}_{e^{i\alpha_{2}}}(h_2,h_2')=L^{1}_{e^{-i\alpha_{2}}}(h_2,h_2')=\mu_{r}.
            \end{equation*}
            By Statement (1) of Lemma \ref{lemma for the injectivity of the landslide}, we obtain $h_1=h_1'$ and $h_2=h_2'$. By Statement (2) of Lemma \ref{lemma for the injectivity of the landslide}, we get that
            \begin{equation*}
                h_1=L^{1}_{e^{i\alpha_{1}}}(h_1,h_1')=\mu_{l}=L^{1}_{e^{-i\alpha_{1}}}(h_1,h_1')=h_2.
            \end{equation*}
            This implies that $(h_1,h_1')=(h_2,h_2')$. The proof is complete.
        \end{proof}

        \begin{remark}
            Note that Theorem \ref{tm:K1K2} also holds for the non-singular case. This implies that for a non-Fuchsian convex GHM AdS manifold $N$ (with particles or not), any two spacelike surfaces of distinct constant curvatures are not isotopic.
        \end{remark}



\end{document}